\documentclass[a4paper,10pt,reqno]{amsart}
\usepackage{amsmath,amsfonts,amsthm}
\usepackage[mathscr]{eucal}
\usepackage{rotating}
\usepackage{multirow}
\usepackage{slashbox}

\numberwithin{equation}{section}
\newtheorem{theorem}{Theorem}[section]
\newtheorem{lemma}[theorem]{Lemma}
\newtheorem{remark}[theorem]{Remark}
\newtheorem{proposition}[theorem]{Proposition}

\theoremstyle{definition}
\newtheorem{definition}[theorem]{Definition}
\newtheorem{assumption}[theorem]{Assumption}
\newtheorem{example}[theorem]{Example}

\begin{document}

\title[Existence of Neutral Equations]{On the Existence and Uniqueness of Solutions of Stochastic Equations of Neutral Type}

\author{John A. D. Appleby}
\address{Edgeworth Centre for Financial Mathematics, School of Mathematical Sciences, Dublin City University,
Ireland} \email{john.appleby@dcu.ie} \urladdr{webpages.dcu.ie/\textasciitilde applebyj}
\author{Huizhong Appleby--Wu}
\address{Department of Mathematics, St Patrick's College, Drumcondra, Dublin 9, Ireland} \email{huizhong.applebywu@spd.dcu.ie}
\author{Xuerong Mao}
\address{Department of Mathematics and Statistics, University of Strathclyde,
Livingstone Tower, 26 Richmond Street, Glasgow G1 1XT, United
Kingdom.}
\email{x.mao@strath.ac.uk}
\urladdr{http://personal.strath.ac.uk/x.mao/}
%\author{Huizhong Wu}
%\address{School of Mathematical Sciences, Glasnevin, Dublin 9, Dublin City University,
%Ireland} \email{huizhong.wu4@mail.dcu.ie}

%\thanks{The first author was partially
%supported by an Albert College Fellowship awarded by Dublin City
%University's Research Advisory Panel. The second author was funded by Science Foundation Ireland (SFI) under the
%Research Frontiers Programme grant RFP/MAT/0018 ``Stochastic
%Functional Differential Equations with Long Memory''.}

\thanks{We gratefully acknowledge the support of this work by Science Foundation Ireland (SFI) under the
Research Frontiers Programme grant RFP/MAT/0018 ``Stochastic
Functional Differential Equations with Long Memory''. JA also thanks
SFI for the support of this research under the Mathematics
Initiative 2007 grant 07/MI/008 ``Edgeworth Centre for Financial
Mathematics''.}



\date{2 May 2013}

\begin{abstract}
This paper considers some the existence and uniqueness of strong
solutions of stochastic neutral functional differential equations.
The conditions on the neutral functional relax those commonly used
to establish the existence and uniqueness of solutions of NSFDEs for
many important classes of functional, and parallel the conditions
used to ensure existence and uniqueness of solutions of related
deterministic neutral equations. Exponential estimates on the almost
sure and p--th mean rate of growth of solutions under the weaker
existence conditions are also given.
\end{abstract}

\maketitle

\section{Introduction}
Over the last ten years, a body of work has emerged concerning the
properties of stochastic neutral equations of It\^o type. Of course,
one of the most fundamental questions is whether solutions of such
equations exist and are unique. A great many of these results have
been established by Mao and co-workers.

In this paper, we concentrate for simplicity on autonomous
stochastic neutral functional differential equations, and establish
existence and uniqueness of solutions under weaker conditions than
currently extant in the literature. The solutions will be unique
within the class of continuous adapted processes, and will also
exist on $[0,\infty)$. Also for simplicity, we assume that all
functionals are globally linearly bounded and globally Lipschitz
continuous (with respect to the sup--norm topology). The most
general finite--dimensional neutral equation of this type is
\begin{gather}\label{eq:main}
d(X(t)-D(X_t))=f(X_t)\,dt + g(X_t)\,dB(t), \quad 0\leq t\leq T;\\
X(t)=\psi(t),\quad t\in [-\tau,0].
\end{gather}
where $\tau>0$, $\psi\in C([-\tau,0];\mathbb{R}^d)$, $B$ is an
$m$--dimensional standard Brownian motion, $D$ and $f$ are
functionals from $C([-\tau,0];\mathbb{R}^d)$ to $\mathbb{R}^d$ and
$g:C([-\tau,0];\mathbb{R}^d\times \mathbb{R}^m)\to \mathbb{R}^d$. It
is our belief that the results presented in this paper can be
extended to non--autonomous equations, to equations which obey only
local Lipschitz continuity conditions, and to equations with local
linear growth bounds. Naturally, in these circumstances, we cannot
expect solutions to necessarily be global; instead, one can talk
only about the existence of local solutions.

To the best of the authors' knowledge, all existing existence
results concerning stochastic neutral equations in general, and
\eqref{eq:main} in particular, involve a ``contraction condition''
on the operator $D$ on the righthand side. We term the operator $D$
the \emph{neutral functional} throughout this paper, and the
functional $E:C([-\tau,0];\mathbb{R}^d)\to \mathbb{R}^d$ defined by
$E(\phi):=\phi(0)-D(\phi)$ the \emph{neutral term}. The
contraction condition on $D$ is that there exists a $\kappa\in
(0,1)$ such that
\begin{equation} \label{eq.contrcondnmao}
|D(\phi)-D(\varphi)|\leq \kappa\|\phi-\varphi\|_{\text{sup}},
\quad\text{ for all $\phi$, $\varphi\in C([-\tau,0];\mathbb{R}^d)$},
\end{equation}
where $|\cdot|$ is the standard norm in $\mathbb{R}^d$ and
$\|\phi\|_{\text{sup}}:=\sup_{-\tau\leq s\leq 0}|\phi(s)|$ where
$\phi\in C([-\tau,0];\mathbb{R}^d)$. Under this condition, as well
as conventional Lipschitz conditions on $f$ and $g$, it can be shown
that \eqref{eq:main} has a unique continuous adapted solution on
$[0,T]$ for every $T>0$.

While the condition \eqref{eq.contrcondnmao} is certainly sufficient
to ensure existence and uniqueness of solutions, until now it has
not been understood whether this condition is necessary. However,
comparison with the existence theory for the deterministic neutral
equation corresponding to \eqref{eq:main} viz.,
\begin{gather}\label{eq:maindet}
\frac{d}{dt}(x(t)-D(x_t))=f(x_t), \quad 0\leq t\leq T;\\
x(t)=\psi(t),\quad t\in [-\tau,0].
\end{gather}
would lead one to suspect that the condition
\eqref{eq.contrcondnmao} is too strong, at least in some
circumstances. To take a simple scalar example, suppose that
$f:C([-\tau,0];\mathbb{R})\to\mathbb{R}$ is globally Lipschitz
continuous, and that $w\in C([-\tau,0];\mathbb{R}^+)$ is such that
\begin{equation}  \label{eq.introweightgt1}
\int_{-\tau}^0 w(s)\,ds >1.
\end{equation}
Then the solution of
\begin{gather*}
\frac{d}{dt}(x(t)-\int_{-\tau}^0 w(s)x(t+s)\,ds=f(x_t), \quad 0\leq t\leq T;\\
x(t)=\psi(t),\quad t\in [-\tau,0].
\end{gather*} exists and is
unique in the class of continuous functions. On the other hand,
extant results do not enable us to make a definite conclusion
concerning the existence and uniqueness of solutions of
\begin{gather}\label{eq.stochexampintro}
d(X(t)-\int_{-\tau}^0 w(s)X(t+s)\,ds)=f(X_t)\,dt + g(X_t)\,dB(t), \quad 0\leq t\leq T;\\
X(t)=\psi(t),\quad t\in [-\tau,0].
\end{gather}
when $g:C([-\tau,0];\mathbb{R})\to\mathbb{R}$ is also globally
Lipschitz continuous, because the functional $D$ defined by
\begin{equation}\label{eq.introroguefunct}
D(\phi)=\int_{-\tau}^0 w(s)\phi(s)\,ds
\end{equation}
does not obey \eqref{eq.contrcondnmao} if $w$ obeys \eqref{eq.introweightgt1}.

It transpires that the condition of uniform non--atomicity at zero of the functional $D$, which was introduced by
Hale in the deterministic theory, and ensures the existence and uniqueness of a solution of the equation \eqref{eq:maindet},
also ensures the existence of a unique solution of \eqref{eq:main}, under Lipschitz continuity conditions on $f$ and $g$. We
discuss this non--atomicity condition presently, but note that it entails the existence of a number $s_0\in(0,\tau)$ and a non--decreasing function
$\kappa:[0,s_0]\to \mathbb{R}$ such that $\kappa(s_0)<1$ and
\begin{multline} \label{eq.contrcondnhale}
|D(\phi)-D(\varphi)|\leq \kappa(s)\|\phi-\varphi\|_{\text{sup}}
\quad\text{ for all $\phi$, $\varphi\in C([-\tau,0];\mathbb{R}^d)$},\\
\text{such that $\phi=\varphi$ on $[-\tau,-s]$ and $s\in[0,s_0]$.}
\end{multline}
Roughly speaking, it can be seen that \eqref{eq.contrcondnhale} relaxes \eqref{eq.contrcondnmao} by allowing the functions $\phi$ and $\varphi$
to be equal on a subinterval of $[-\tau,0]$, thereby effectively reducing the Lipschitz constant in \eqref{eq.contrcondnmao} from a number
greater than unity to a number less than unity. As an example, the functional in \eqref{eq.introroguefunct} obeys \eqref{eq.contrcondnhale}
even under the condition \eqref{eq.introweightgt1} on $w$. Therefore, we can conclude that \eqref{eq.stochexampintro} has
a unique solution; existing results would however require $w$ to obey $\int_{-\tau}^0 w(s)\,ds<1$.

The condition \eqref{eq.contrcondnmao} has to date played a very
important role in the analysis of properties of solutions of
\eqref{eq:main}. It is a key assumption in proofs of estimates on
the almost sure and $p$-th mean rate of growth of solutions of
\eqref{eq:main}. It is also required in results which deal with the
almost sure and $p$--th mean asymptotic stability of solutions.
Results on the $L^p$ continuity of solutions, and even results on
numerical methods to simulate the solution of \eqref{eq:main}, rely
on the condition \eqref{eq.contrcondnmao}. However, corresponding
results for the underlying deterministic equation \eqref{eq:maindet}
regarding asymptotic behaviour, regularity of solutions, and
numerical methods can be established under the weaker condition
\eqref{eq.contrcondnhale}.

It is therefore reasonable to ask whether fundamental results on
e.g., asymptotic behaviour, can still be established for solutions
of \eqref{eq:main} under the weaker condition
\eqref{eq.contrcondnhale}, which is shown in this paper to be
sufficient to ensure solutions exist. Towards this end, in this
paper we prove results on almost sure and $p$--th mean exponential
estimates on the growth of the solution of \eqref{eq:main} using the
condition \eqref{eq.contrcondnhale} in place of
\eqref{eq.contrcondnmao}. Although we confine our attention here to
the study of these exponential estimates, it is of obvious interest
to investigate further the properties of solutions of stochastic
neutral equations under the weaker non--atomicity condition
\eqref{eq.contrcondnhale} which have, owing to the absence of
existence results, remained unconsidered until now.

%without placing a ``contraction condition'' on the neutral operator
%$D$. Such a condition takes the form  To date, existence and
%uniqueness results for general SNFDEs have relied on such a
%condition, in contrast to the situation for deterministic neutral
%functional differential equations. The results are stated and proved
%for autonomous equations but can readily be extended to classes of
%non--autonomous equations. Moreover, we have imposed global
%Lipschitz conditions, but these can easily be supplanted by local
%conditions if we require only local existence and uniqueness of
%solutions.
%
%\section{Introduction}
%This paper determines conditions for the existence and uniqueness of strong solutions of stochastic neutral functional differential equation
%\begin{gather}\label{eq:main}
%d(X(t)-D(X_t))=f(X_t)\,dt + g(X_t)\,dB(t), \quad 0\leq t\leq T;\\
%X(t)=\psi(t),\quad t\in [-\tau,0].
%\end{gather}
%without placing a ``contraction condition'' on the neutral operator $D$. Such a condition takes the form  To date, existence and uniqueness results for general SNFDEs have relied
%on such a condition, in contrast to the situation for deterministic neutral functional differential equations. The results are stated and proved for autonomous equations but can readily be extended to classes of non--autonomous equations. Moreover, we have imposed global Lipschitz conditions, but
%these can easily be supplanted by local conditions if we require only local existence and uniqueness of solutions.
%
%The removal of the ``contraction condition'' is predicated on the neutral functional being of the form

Neutral delay differential equations have been used to describe
various processes in physics and engineering sciences \cite{HaleLun93}, \cite{Stepan:1989}. For
example, transmission lines involving nonlinear boundary conditions
\cite{Hale:77}, cell growth dynamics \cite{BakBochPaul:1998}, propagating pulses in cardiac tissue
\cite{CourGlassKeen:1993} and drillstring vibrations \cite{BalJansMcClint:2003} have been described by means of
neutral delay differential equations.

\textbf{To do}
\begin{itemize}
\item Almost sure exponential growth bound.
\end{itemize}

\section{Mathematical Preliminaries}

In this section, we introduce some notation that will be used throughout the paper, state and comment on known results on the existence of solutions
of the stochastic neutral equation \eqref{eq:main}, and introduce in precise terms the weaker conditions used here on the neutral functional $D$ which will still guarantee existence and uniqueness of solutions of \eqref{eq:main}.

\subsection{Notation}
We denote the upper Dini derivative by $D^+$, i.e. if $f:\mathbb{R}\to\mathbb{R}$ is
continuous, then
\[
D^+ f(t):=\limsup_{h\to 0^+}\frac{f(t+h)-f(t)}{h}.
\]
For any  $d\in\mathbb{N}$ and $\tau>0$, we define $C([-\tau, 0];\mathbb{R}^d)$ is the space of continuous functions from $[-\tau, 0] \to \mathbb{R}^d$ with sup norm. The sup norm $\|\cdot\|_\text{sup}$ on $C([-\tau,0];\mathbb{R}^d)$ is defined so that for $\phi\in C([-\tau,0];\mathbb{R}^d)$ we have
\[
\|\phi\|_\text{sup}=\max_{-\tau\leq s\leq 0} |\phi(s)|,
\]
where $|\cdot|$ denotes the usual Euclidean norm on $\mathbb{R}^d$.

Let $\phi$ be a function from $[-\tau, t_1] \to \mathbb{R}^d$. Let $t\in [0, t_1] \subset \mathbb{R}$. We use $\phi_t$ to denote the function on $[-\tau, 0]$ defined by $\phi_t(s) = \phi(t + s)$ for $-\tau \leq s \leq 0$.

Let $d,d^\prime $ be some positive integers and $\mathbb{R}^{d\times
d^\prime}$ denote the space of all $d \times d^\prime$ matrices
with real entries.
%The identity matrix on $\mathbb{R}^{d\times d}$ is denoted by $\Id_d$.
We equip $\mathbb{R}^{d\times d^\prime}$ with a norm
$|\cdot|$ and write $\mathbb{R}^d$ if $d^\prime=1$ and $\mathbb{R}$ if
$d=d^\prime=1$. We denote by $\mathbb{R}^+$ the half-line $[0,\infty)$.
%The complex plane is denoted by $\C$ and $\C_0:=\{z\in\C:\, \Re z\ge 0\}$.

Let $M(\mathbb{R}^+,\mathbb{R}^{d\times d^\prime})$ be the space of finite Borel
measures on $\mathbb{R}^+$ with values in $\mathbb{R}^{d\times d^\prime}$. The
total variation of a measure $\nu$ in $M(\mathbb{R},\mathbb{R}^{d\times
d^\prime})$ on a Borel set $B\subseteq \mathbb{R}^+$ is defined by
\begin{align*}
 |\nu|\!(B):=\sup\sum_{i=1}^N |\nu(E_i)|,
\end{align*}
where $(E_i)_{i=1}^N$ is a partition of $B$ and the supremum is
taken over all partitions. The total variation defines a positive
scalar measure $|\nu|$ in $M(\mathbb{R}^+,\mathbb{R})$. If one specifies
temporarily the norm $|\cdot|$ as the $l^1$-norm on the space
of real-valued sequences and identifies $\mathbb{R}^{d\times d^\prime}$ by
$\mathbb{R}^{dd^\prime}$ one can easily establish for the measure
$\nu=(\nu_{i,j})_{i,j=1}^d$ the inequality
\begin{align}\label{eq.totalvarest}
 |\nu|\!(B)\leq C \sum_{i=1}^d\sum_{j=1}^d |\nu_{i,j}|\!(B)
 \qquad\text{for every Borel set }B\subseteq \mathbb{R}^+
\end{align}
with $C=1$. Then, by the equivalence of every norm on
finite-dimensional spaces, the inequality \eqref{eq.totalvarest}
holds true for the arbitrary norms $|\cdot|$ and some constant
$C>0$. Moreover, as in the scalar case we have the fundamental
estimate
\begin{align*}
 \left|\int_{\mathbb{R}^+} \nu(ds)\, f(s)\right| \leq \int_{\mathbb{R}^+} |f(s)|\,|\nu|\!(du)
\end{align*}
for every function $f:\mathbb{R}^+\to\mathbb{R}^{d^\prime \times d^{\prime\prime}}$
which is $|\nu|$-integrable. The convolution of a function $f$
and a measure $\nu$ is defined by
\begin{align*}
 \nu\ast f:\mathbb{R}^+\to \mathbb{R}^{d\times d^{\prime\prime}},\qquad
  (\nu\ast f)(t):=\int_{[0,t]} \nu(ds)\, f(t-s).
\end{align*}
The convolution of two functions is defined analogously.

\subsection{Existing Results for Stochastic Neutral Equations}
Let $m$ and $d$ be positive integers. Let $(\Omega,\mathcal{F},\mathbb{P})$ be a complete probability space with the filtration $(\mathcal{F}(t))_{t\geq 0}$ satisfying the usual conditions.

Let $B=\{B(t):t\geq 0\}$ be an $m$--dimensional Brownian motion defined on the space.
Let $\tau>0$ and $0<T<\infty$. Let the functionals $D$, $f$ and $g$ defined by
\[
D:C([-\tau,0];\mathbb{R}^d)\to \mathbb{R}^d, \quad
f:C([-\tau,0];\mathbb{R}^d)\to \mathbb{R}^d, \quad
g:C([-\tau,0];\mathbb{R}^d)\to \mathbb{R}^{d\times m}
\]
all be Borel--measurable.

Consider the $d$--dimensional neutral stochastic functional differential equation
\begin{equation} \label{eq.maineqdiff}
d(X(t)-D(X_t))=f(X_t)\,dt + g(X_t)\,dB(t), \quad 0\leq t\leq T.
\end{equation}
This should be interpreted as the integral equation
\begin{equation} \label{eq.maineqint}
X(t)-D(X_t)=X(0)-D(X_0)+\int_0^t f(X_s)\,ds + \int_0^t g(X_s)\,dB(s), \quad \text{for all $t\in[0,T]$}.
\end{equation}
For the initial value problem we must specify the initial data on the interval $[-\tau,0]$ and hence we impose the initial
condition
\begin{equation} \label{eq.ic}
X_0=\psi=\{\psi(\theta):-\tau\leq \theta\leq 0\}\in L^2_{\mathcal{F}(0)}([-\tau,0];\mathbb{R}^d),
\end{equation}
that is $\psi$ is an $\mathcal{F}(0)$--measurable $C([-\tau,0];\mathbb{R}^d)$--valued random variable such that $\mathbb{E}[|\psi|^2]<+\infty$.
The initial value problem for equation \eqref{eq.maineqdiff} is to find the solution of \eqref{eq.maineqdiff} satisfying the initial data
\eqref{eq.ic}. We give the definition of the solution in this context
\begin{definition}
An $\mathbb{R}^d$--valued stochastic process $X=\{X(t):-\tau\leq t\leq T\}$ is called a solution to equation \eqref{eq.maineqdiff}
with initial data \eqref{eq.ic} if it has the following properties:
\begin{itemize}
\item[(i)] $t\mapsto X(t,\omega)$ is continuous for almost all $\omega\in \Omega$ and $X$ is $(\mathcal{F}(t))_{t\geq 0}$--adapted;
\item[(ii)] $\{f(X_t)\}\in L^1([0,T];\mathbb{R}^d)$ and $\{g(X_t)\}\in L^2([0,T];\mathbb{R}^{d\times m})$;
\item[(iii)] $X_0=\psi$ and \eqref{eq.maineqint} holds for every $t\in[0,T]$.
\end{itemize}
A solution $X$ is said to be unique if any other solution $\bar{X}$ is indistinguishable from it i.e.,
\[
\mathbb{P}[X(t)=\bar{X}(t) \text{ for all $-\tau\leq t\leq T$}]=1.
\]
\end{definition}
We now make the following assumptions on the functionals $f$ and $g$ in order to ensure the existence and uniqueness of
solutions of \eqref{eq.maineqdiff}. They will hold throughout the paper.
\begin{assumption} \label{ass.fgliplinbdd}
There exists $K>0$ such that for all $\phi, \varphi\in C([-\tau,0];\mathbb{R}^d)$
\begin{equation} \label{eq.fglip}
|f(\varphi)-f(\phi)|\leq K\|\varphi-\phi\|_\text{sup}, \quad ||g(\varphi)-g(\phi)||\leq K\|\varphi-\phi\|_\text{sup}.
\end{equation}
There exists $\bar{K}>0$ such that for all $\phi, \varphi\in C([-\tau,0];\mathbb{R}^d)$
\begin{equation} \label{eq.fglinboud}
|f(\varphi)|\leq \bar{K}(1+\|\varphi\|_\text{sup}), \quad ||g(\varphi)||\leq \bar{K}(1+\|\varphi\|_\text{sup}).
\end{equation}
\end{assumption}
The following result is Theorem 6.2.2 in \cite{Mao97}; it concerns the existence and uniqueness of solutions of the stochastic neutral functional differential equation \eqref{eq.maineqdiff}.
\begin{theorem} \label{thm.maoexistmain}
Suppose that the functionals $f$ and $g$ obey \eqref{eq.fglip} and \eqref{eq.fglinboud} and that the functional $D$
obeys
\begin{multline} \label{eq.dlipcontract}
\text{There exists $\kappa\in(0,1)$ such that for all $\phi, \varphi\in C([-\tau,0];\mathbb{R}^d)$}\\
|D(\varphi)-D(\phi)|\leq \kappa\|\varphi-\phi\|_{\text{sup}}.
\end{multline}
Then there exists a unique solution $X$ to \eqref{eq.maineqdiff} with initial data \eqref{eq.ic}.
Moreover the solution belongs to $\mathcal{M}^2([-\tau,T];\mathbb{R}^d)$.
\end{theorem}
On the other hand, a restriction of this type on the neutral functional $D$ such as \eqref{eq.dlipcontract} is not needed in the case when it depends purely on delayed arguments. See \cite[Theorem 6.3.1]{Mao97}.

\subsection{Assumptions on the Neutral Functional}
In order to orient the reader to the question of existence which is addressed in this paper, we must first introduce some results and
notation from the theory of deterministic neutral differential equations.
Consider systems of nonlinear functional differential equations
of neutral type having the form
\begin{equation} \label{eq.neutral chukwu}
\frac{d}{dt} E(x_t) = f(x_t),
\end{equation}
where the operator $E : C \to \mathbb{R}^d$ is \emph{atomic} at $0$ and \emph{uniformly atomic} at
$0$ in the sense of Hale~\cite[pp 170--173]{Hale:71}, and where $f :C\to \mathbb{R}^d$ is
continuous and uniformly Lipschitzian in the last argument. In \eqref{eq.neutral chukwu}, instead of the atomicity assumption on $E$, we
may assume that $E$ is of the form
\[
E(\phi)=\phi(0)-D(\phi)
\]
where $D : C\to \mathbb{R}^d$ is continuous and is uniformly nonatomic at zero on $C$
in the following sense.

\begin{definition} \label{def.unifnonatomic}
For any $\phi \in C$, and
$s \geq 0$, let
\[
Q(\phi, s) = \{\varphi \in  C :  \varphi(\theta) = \phi(\theta),
\theta < -s, \theta \in [-\tau,0]\}.
\]
We say that a continuous function $D: C \to \mathbb{R}^d$ is \emph{uniformly nonatomic at zero} on $C$ if, for any $\phi \in C$, there exist $T_1 > 0$, independent of $\phi$, and a positive scalar function $\rho(\phi, s)$, defined for $\phi \in C$, $0 \leq  s \leq  T_1$, nondecreasing in $s$ such that
\begin{equation} \label{eq.rho0rho0lt1}
\rho_0(s):=
\sup_{\phi\in C} \rho(\phi, s), \quad \rho_0(T_1)=:k<1,
\end{equation}
and
\begin{equation} \label{eq.nonatomlip2a}
|D(\varphi_1)-D(\varphi_2)|\leq \rho_0(s) \|\varphi_1-\varphi_2\|_\text{sup}, \text{ for $\varphi_1,\,\varphi_2 \in Q(\phi, s)$ and all $0 \leq s \leq T_1$}.
\end{equation}
\end{definition}
We note that the definition implies both that $\rho_0$ is non--decreasing and that $\rho_0$ is independent of $\phi$. Therefore
a consequence of \eqref{eq.nonatomlip2a} is
\begin{multline} \label{eq.nonatomlip2}
|D(\varphi_1)-D(\varphi_2)|\leq \rho_0(s) \|\varphi_1-\varphi_2\|_\text{sup}, \text{ for $\varphi_1,\,\varphi_2 \in Q(\phi, s)$},\\
 \text{ and all $0 \leq s \leq T_1$ and all $\phi\in C$}.
\end{multline}
We tend to use this consequence of the definition in practice.

It is instructive to compare the conditions \eqref{eq.rho0rho0lt1} and \eqref{eq.nonatomlip2a} with Mao's condition \eqref{eq.dlipcontract}
on the neutral functional $D$. We first note that \eqref{eq.dlipcontract} implies both \eqref{eq.rho0rho0lt1} and \eqref{eq.nonatomlip2a}
and so implies that $D$ is uniformly nonatomic at 0 in $C([-\tau,0];\mathbb{R}^d)$, so that \eqref{eq.dlipcontract} is not a weaker condition
that uniform nonatomicity. Indeed, as shown by the functional given in \eqref{eq.introroguefunct}, the condition \eqref{eq.dlipcontract} is
a strictly stronger condition.

It is known (\cite{ChukwuSimp:1989,Hale:71,HaleCruz:70}) that under these assumptions on $D$, and $f$ for
each $\phi\in C$ there is a unique solution of \eqref{eq.neutral chukwu} with initial value $\phi$ at $0$. The solution is continuous with respect to initial data. For definition of solutions see~\cite{HaleCruz:70}. In the sequel $t_1$ is fixed and is in the interval of definition $[0, T]$, of solutions of \eqref{eq.neutral chukwu}.

We make the following related assumption on the functional.

\begin{assumption}\label{assD}
Let $\phi\in C([-\tau,0];\mathbb{R}^d)$ and assume $D_0$, $D_1 : C\to \mathbb{R}^d$ such that
\begin{equation}\label{Ddecomp}
D(\phi)=D_0(\phi)+D_1(\phi).
\end{equation}
Suppose there exists $\delta>0$ and $H: C([-\tau, 0]; \mathbb{R}^d)\to \mathbb{R}^d$ such that
\begin{equation}\label{D0}
D_0(\phi):=H(\{\phi(s):\,-\tau\leq s\leq -\delta<0\}), \\ \nonumber \text{for all $\phi\in C([-\tau, 0];\mathbb{R}^d)$}.
\end{equation}
Suppose further that $D_1$ is uniformly non--atomic at zero on $C$, so that there exists $0<T_1\leq \delta$ and $k\in(0,1)$ as given
in definition \ref{def.unifnonatomic} such that \eqref{eq.rho0rho0lt1} and \eqref{eq.nonatomlip2} hold.
%\begin{gather}\label{D1}
%\forall \,t\in[nT_1,(n+1)T_1], \quad n \in\mathbb{N}\cup \{0\},\quad x,y \in C([-\tau, T]; \mathbb{R}^d), \\
%\nonumber\text{and}\quad x(s)=y(s), \quad -\tau\leq s \leq nT_1,\\
%\nonumber |D_1(x_t)-D_1(y_t)|\leq k \sup_{-T_1\leq s \leq 0}|x(t+s)-y(t+s)|.
%\end{gather}
\end{assumption}
We can choose $T_1<\delta$ without loss of generality in order to ensure that the pure delay functional $D_0$ which depends on $\phi\in C([-\tau,0];\mathbb{R}^d)$ only on $[-\tau,-\delta]$ does not interact with the functional $D_1$ which can depend on $\phi$ on all $[-\tau,0]$. One consequence of the decomposition
of $D$ in \eqref{Ddecomp} is that the continuity condition on $k$ required in Hale's definition of uniform non--atomicity can be dropped.

We make a linear growth assumption on $D$ which is slightly non--standard also.
\begin{assumption}\label{assDlineargr}
%For all $x\in C([-\tau, T];\mathbb{R}^d)$, there exist $k\in(0,1)$ and $K_D>0$ such that
For all $\phi\in C([-\tau, 0];\mathbb{R}^d)$, there exist $k\in(0,1)$ and $K_D>0$ such that
\begin{gather}
%\label{1stint} \forall\, t\in[0,T_1],\quad |D(x_t)|\leq K_D+k\sup_{-\tau\leq s\leq 0}|x(t+s)|;\\
%\label{nint} \forall\,t\in[nT_1, (n+1)T_1],\quad n\in \mathbb{N}, \\ \nonumber
%|D(\phi_t)|\leq K_D(1+\sup_{-\tau\leq s\leq -T_1}|\phi(t+s)|)+k\sup_{-T_1\leq s\leq 0}|\phi(t+s)|.
\label{nint}
|D(\phi)|\leq K_D(1+\sup_{-\tau\leq s\leq -T_1}|\phi(s)|)+k\sup_{-T_1\leq s\leq 0}|\phi(s)|.
 %\quad \text{for all $\phi\in C([-\tau,0];\mathbb{R}^d)$}.
\end{gather}
\end{assumption}
The numbers $k$ and $T_1$ can be chosen to be the same as those in Assumption \ref{assD} without loss of generality, and we choose to do
so. One reason for this is that the choice that $T_1<\delta$ in Assumption \ref{assD} ensures that the pure delay functional $D_0$ does not make a contribution to the constant
$k$ in the second term on the right hand side of \eqref{nint} which might force $k>1$.
The linear growth bound on $D(\phi)$ arising from the dependence on $\phi$ over the interval $[-\tau,-T_1]$ guarantees the existence of second moments of the solution of \eqref{eq:main}. Notice that no restriction is made on the size of the constant $K_D$, while we require $k\in(0,1)$.

%%%%%%%%%%%%%%%%%%%%%%%%%%%%%%%%%%%%%%%%%%%%%%%%%%%%%%%%%%%%%%%%%%%%%%%%%%%%%%%%%%%%%%%%%%%%%%%%%%%%%%%%%%%%%%%%%%%%%%%%%%%%%%%%%%%%%%%%%%%%%%%%%%%%
\section{Discussion of Main Results}\label{sec:existence}
In this section we state and discuss the main results of the paper. We state our main existence result, and give examples of functionals
to which it applies. We then show, under the condition that $D$ is uniformly non--atomic at zero in $C([-\tau,0];\mathbb{R}^d)$, that the solution
$X$ of \eqref{eq.maineqdiff} enjoys exponential growth bounds in both a $p$--th mean and almost sure sense. Finally, we give examples of
equations for which the neutral functional $D$ is not uniformly non--atomic at zero, and for which solutions of \eqref{eq.maineqdiff}
do not exist.
\subsection{Existence result}
The main result of this paper relaxes the contraction constant in \eqref{eq.dlipcontract} in the case when the functional $D$ is
composed of a mixture of pure delay and instantaneously interacted functional. For any $T>0$ and $\tau\geq 0$ we define 
$\mathcal{M}^2([-\tau,T];\mathbb{R}^d)$ to be the space of all $\mathbb{R}^d$--valued adapted process $U=\{U(t):-\tau\leq t\leq T\}$ such that 
\[
\mathbb{E}\left[\sup_{-\tau\leq s\leq T}|U(s)|^2\right] < +\infty.
\]
\begin{theorem}\label{thmexiuni}
Suppose that the functionals $D$ obeys Assumption \ref{assD} and Assumption \ref{nint}, $f$ and $g$ obey
Assumption \ref{ass.fgliplinbdd}. Then there exists a unique solution to equation \eqref{eq:main}.
Moreover, the solution is in $\mathcal{M}^2([-\tau,T];\mathbb{R}^d)$.
\end{theorem}
We now give two examples to which Theorem \ref{thmexiuni} can be applied.
\begin{example} \label{examp:neut1}
Consider the neutral functional $D$ defined by
\begin{equation} \label{eq.uniffunctionalgenex}
D(\varphi)=h_0(\varphi(0))+\sum_{i=1}^N h_i(\varphi(-\tau_i))+\int_{[-\tau_0,0]}w(s)h(\varphi(s))\,ds,
\end{equation}
where $\varphi \in C([-\max_{i\geq 1}\{\tau_i\}\vee \tau_0, 0];\mathbb{R}^d)$; $h$ is global Lipschitz continuous and linearly bounded; $w$ is continuous; For each $i\in\mathbb{N}$, $\tau_i>0$, $h_i$ is continuous and globally linearly bounded. It is easy to see that under either of the following two conditions, a unique solution exists:
\begin{itemize}
\item [(i)] If $h_0$ is also global Lipschitz continuous and linearly bounded, moreover, for any $x, y \in\mathbb{R}^d$, $|h_0(x)-h_0(y)|\leq k |x-y|$ with $0<k<1$.
\item [(ii)] If $h_0(x)=Ax$, $A\in \mathbb{R}^{d\times d}$ and $\text{det}(I-A)\neq 0$. In this case, equation \eqref{eq:main} can be rearranged by dividing both sides by $(I-A)^{-1}$ to obtain a unique solution regardless the value of $k$.
\end{itemize}
The two cases illustrate the importance of both invertibility and non-atomicity in ensuring a unique solution of equation \eqref{eq:main}.
\end{example}
\begin{example} \label{examp:neut2}
Consider $D(\varphi)=K\max_{-\tau\leq s\leq -\tau'}||\varphi(s)||$ where $0\leq \tau'<\tau$. If $\tau'>0$, then for all $K\in\mathbb{R}$, a unique solution exists. In this case, $D$ plays the role of $D_0$ in \eqref{eq:main}. However, if $\tau'=0$, then we require that $|K|<1$.
\end{example}
\subsection{Exponential estimates on the solution}
In this subsection we state our results on the existence of moment and almost sure exponential estimates on the solution of \eqref{eq:main}.
Results of this kind have been proven by Mao in~\cite[Chapter 6]{Mao2008} under the condition \eqref{eq.dlipcontract}. However, in
this paper we establish similar estimates under the weaker assumption that $D$ is uniformly non--atomic at zero. In our proof, this relaxation of
the condition comes at the expense of a strengthening of our hypotheses on the functionals $D$, $f$ and $g$. The new hypotheses, which tend to preclude the functionals being closely related to maximum functionals, are nonetheless very natural for equations with point or distributed delay.
The proofs rely on differential and integral inequalities, in contrast to those in~\cite[Chapter 6]{Mao2008}.
\begin{theorem}\label{thmexpest}
Suppose that $f$ and $g$ are globally Lipschitz continuous and that $D$ is uniformly non--atomic at zero. Then there exists
a unique continuous solution $X$ of equation \eqref{eq:main}. Suppose further that there exist positive real numbers $C_f$, $C_g$ and $C_D$ such that
\begin{gather}
\label{f1}|f(\varphi)|\leq C_f+\int_{[-\tau,0]}\nu(ds)|\varphi(s)|;\\
\label{g1}||g(\varphi)||\leq C_g+\int_{[-\tau,0]}\eta(ds)|\varphi(s)|;\\
\label{D1}|D(\varphi)|\leq C_D+\int_{[-\tau,0]}\mu(ds)|\varphi(s)|,
\end{gather}
where  $\nu$, $\eta$ and $\mu$ $\in M([-\tau,0];\mathbb{R}^+)$. Let $p\geq 2$, $\varepsilon>0$ and define
\[
\beta_1=\beta_1(p,\varepsilon):=\frac{\varepsilon p(p-1)}{2},\quad \lambda(du)=\lambda_{p,\varepsilon}(du):=\nu(du)\cdot\frac{1}{\varepsilon^{p-1}}+\eta(du)\cdot \frac{p-1}{\varepsilon^{\frac{p-2}{2}}}.
\]
Then there exists a positive real number $\delta=\delta(p,\varepsilon)$ such that $X$ obeys
\begin{equation}\label{expXp}
\limsup_{t\to\infty}\frac{1}{t}\log{\mathbb{E}[|X(t)|^p]}\leq \delta+\frac{\varepsilon p(p-1)}{2},
\end{equation}
where $\delta$ satisfies
\[
\int_{[-\tau,0]}e^{(\delta+\beta_1)s}\mu(ds)+\int_0^\tau e^{-\delta s}\int_{[-s,0]}e^{\beta_1 u}\lambda(du)\,ds+\frac{e^{-\delta \tau}}{\delta}\int_{[-\tau,0]}e^{\beta_1 u}\lambda(du)=1.
\]
\end{theorem}
We make no claims about the optimality of the exponent in \eqref{expXp}, although $\varepsilon>0$ could be chosen so as to minimise
$\varepsilon\mapsto \delta(p,\varepsilon)+\beta_1(p,\varepsilon)$ for a given value of $p\geq 2$. In a later work we show that an exact exponent can be determined in the case $p=2$ for a scalar linear stochastic neutral equation.
\begin{remark}
We notice that a functional of a form similar to \eqref{eq.uniffunctionalgenex} satisfies the conditions \eqref{f1}, \eqref{g1} or \eqref{D1}.
Suppose for $i=1,\ldots,N$ that $h_i:\mathbb{R}^{d}\to\mathbb{R}^{d'}$ is globally linearly bounded, and satisfies the bound $|h_i(x)|\leq K_i(1+|x|)$ for $x\in\mathbb{R}^d$, and that $\nu_i\in M([-\tau,0];\mathbb{R}^{d\times d'})$, and let
\[
f(\varphi)=\sum_{i=1}^N \int_{[-\tau_i,0]} \nu_i(ds)h_i(\varphi(s)), \quad \varphi\in C([-\tau,0];\mathbb{R}^d),
\]
where $\tau=\max_{i=1,\ldots,N} \tau_i\in(0,\infty)$.
Then
\begin{align*}
|f(\varphi)|
%&\leq \sum_{i=1}^N \left|\int_{[-\tau_i,0]} \nu_i(ds)h_i(\varphi(s))\right|\\
%&\leq  \sum_{i=1}^N \int_{[-\tau_i,0]} |\nu_i|(ds)|h_i(\varphi(s))|\\
%&\leq  \sum_{i=1}^N \int_{[-\tau_i,0]} K_i|\nu_i|(ds)(1+|\varphi(s))|)\\
&\leq  \sum_{i=1}^N \int_{[-\tau_i,0]} K_i|\nu_i|(ds) + \sum_{i=1}^N \int_{[-\tau_i,0]} K_i|\nu_i|(ds)|\varphi(s)|.
\end{align*}
Now set $C_f=\sum_{i=1}^N \int_{[-\tau_i,0]} K_i|\nu_i|(ds)$ and $\nu(ds):=\sum_{i=1}^N K_i|\nu_i|(ds)$ where we define $\nu_i(E)=0$ for every
Borel set $E\subset [-\tau,-\tau_i)$, so that $f$ obeys \eqref{f1}.
\end{remark}
\begin{remark}
First, we note that the conditions \eqref{f1}, \eqref{g1} and \eqref{D1} imply Assumption \ref{ass.fgliplinbdd} and Assumption \ref{assDlineargr}, with which Lemma \ref{lemafinitemo} can be applied. Second,
for any $p\geq 2$, the conditions \eqref{f1}, \eqref{g1} and \eqref{D1} imply
\begin{gather}
\label{f}|f(\varphi)|^p\leq C_f+\int_{[-\tau,0]}\nu(ds)|\varphi(s)|^p;\\
\label{g}||g(\varphi)||^p\leq C_g+\int_{[-\tau,0]}\eta(ds)|\varphi(s)|^p;\\
\label{D}|D(\varphi)|^p\leq C_D+\int_{[-\tau,0]}\mu(ds)|\varphi(s)|^p,
\end{gather}
respectively for a different set of $C_f$, $C_g$ and $C_D$, and rescaled measures $\nu$, $\eta$ and $\mu$. Therefore, for the reason of convenience, we will be using conditions \eqref{f}, \eqref{g} and \eqref{D} in the proof of Theorem \ref{thmexpest}.
\end{remark}
Theorem~\ref{thmexpest} can be used to prove that $X$ obeys an almost sure exponential growth bound.
\begin{theorem}\label{thmexpestas}
Suppose that $f$ and $g$ are globally Lipschitz continuous and that $D$ is uniformly non--atomic at zero. Then there exists
a unique continuous solution $X$ of equation \eqref{eq:main}. Suppose further that there exist positive real numbers
$C_f$, $C_g$ and $C_D$ such that $f$, $g$ and $D$ obey \eqref{f1}, \eqref{g1} and \eqref{D1} respectively where $\nu$, $\eta$ and $\mu$ $\in M([-\tau,0];\mathbb{R}^+)$. Then there exists $\gamma>0$ such that 
\[
\limsup_{t\to\infty} \frac{1}{t}\log \mathbb{E}[|X(t)|^2] \leq \gamma,
\]
and we have the following estimates for $X$:
\begin{itemize}
\item[(i)] If $\int_{[-\tau,0]} \mu(ds)\geq 1$, then $X$ obeys
\begin{equation}\label{expXas1}
\limsup_{t\to\infty} \frac{1}{t}\log |X(t)|\leq \max(\gamma/2,\theta^\ast), \quad \text{a.s.}
\end{equation}
where $\theta^\ast\geq 0$ is defined by $\int_{[-\tau,0]} e^{\theta^\ast s}\mu(ds)=1$. 
\item[(ii)] If $\int_{[-\tau,0]} \mu(ds)< 1$, then $X$ obeys
\begin{equation} \label{expXas2}
\limsup_{t\to\infty} \frac{1}{t}\log |X(t)|\leq \gamma/2, \quad \text{a.s.}
\end{equation}
\end{itemize}
\end{theorem}

\subsection{Non-existence of Solutions of SNFDEs}
In this section, we give examples of scalar stochastic neutral equation which do not have a solution. To the best of our knowledge, examples of stochastic neutral equations which do not have solutions have not  appeared in the literature to date. Our purpose in constructing such examples is
to demonstrate the importance of the existence conditions \eqref{eq.nonatomlip2} and \eqref{eq.dlipcontract} in ensuring the existence of
solutions. We show that both these sufficient conditions are in some sense sharp in two ways. First, by showing that if either condition
\eqref{eq.nonatomlip2} and \eqref{eq.dlipcontract} is slightly relaxed, then solutions to our examples do not exist. Second, by considering the equations for which solutions do not exist as members of parameterised families of equations, we can show that small changes in the parameters lead to equations which have unique solutions.

We consider both equations with continuously distributed functionals and with maximum type functionals. The first class of equation
shows the condition \eqref{eq.nonatomlip2} cannot readily be improved for equations. The condition \eqref{eq.dlipcontract} is shown
to be quite sharp for equations with max--type functionals.

Let $(\Omega,\mathcal{F},\mathbb{P})$ be a complete probability space with the filtration $(\mathcal{F}(t))_{t\geq 0}$ satisfying the usual conditions. Let $B=\{B(t):t\geq 0\}$ be a one--dimensional Brownian motion defined on the space.
Let $\tau>0$ and $0<T<\infty$.

\subsubsection{Equation with continuously distributed delay}
Let the functional $f$ defined by $f:C([-\tau,0];\mathbb{R})\to \mathbb{R}$ be Borel--measurable.
Let $h\in C(\mathbb{R};\mathbb{R})$, $w\in C^1([-\tau,0];\mathbb{R})$ and $\sigma\neq 0$.
Consider the one--dimensional neutral stochastic functional differential equation
\begin{equation} \label{eq.maineqdiffnonexist}
d\left(\epsilon X(t)+\int_{-\tau}^0 w(s)h(X(t+s))\,ds\right)=f(X_t)\,dt + \sigma \,dB(t), \quad 0\leq t\leq T.
\end{equation}
where $\epsilon\in\mathbb{R}$.
For the initial value problem we must specify the initial data on the interval $[-\tau,0]$ and hence we impose the initial
condition
\begin{equation} \label{eq.icnonexist}
X_0=\psi=\{\psi(\theta):-\tau\leq \theta\leq 0\}\in L^2_{\mathcal{F}(0)}([-\tau,0];\mathbb{R}),
\end{equation}
that is $\psi$ is an $\mathcal{F}(0)$--measurable $C([-\tau,0];\mathbb{R})$--valued random variable such that $\mathbb{E}[|\psi|^2]<+\infty$.
\eqref{eq.maineqdiffnonexist} should be interpreted as the integral equation
\begin{multline} \label{eq.maineqintnonexist}
\epsilon X(t)+\int_{-\tau}^0 w(s)h(X(t+s))\,ds =\int_{-\tau}^0 w(s)h(\psi(s))\,ds
\\+\int_0^t f(X_s)\,ds + \int_0^t \sigma\,dB(s), \quad \text{for all $t\in[0,T]$, a.s.}.
\end{multline}
The initial value problem for equation \eqref{eq.maineqdiffnonexist} is to find the solution of \eqref{eq.maineqdiffnonexist} satisfying the initial data
\eqref{eq.icnonexist}. In this context a solution is an $\mathbb{R}$--valued stochastic process $X=\{X(t):-\tau\leq t\leq T\}$ to equation \eqref{eq.maineqdiffnonexist} with initial data \eqref{eq.icnonexist} if it has the following properties:
\begin{itemize}
\item[(i)] $t\mapsto X(t,\omega)$ is continuous for almost all $\omega\in \Omega$ and $X$ is $(\mathcal{F}(t))_{t\geq 0}$--adapted;
\item[(ii)] $\{f(X_t)\}\in L^1([0,T];\mathbb{R})$;
\item[(iii)] $X_0=\psi$ and \eqref{eq.maineqintnonexist} holds.
\end{itemize}

\begin{proposition} \label{prop.nonexiststochadd}
Let $\tau>0$. Let $h\in C(\mathbb{R};\mathbb{R})$, $w\in
C^1([-\tau,0];\mathbb{R})$, $\psi\in C([-\tau,0];\mathbb{R})$ and
$\sigma\neq 0$. Suppose also that
\[
t\mapsto f(x_t) \text{ is in $C([0,\infty),\mathbb{R})$ for each
$x\in C([-\tau,\infty),\mathbb{R})$}.
\]
Let $T>0$ and $\epsilon=0$. Then there is no process $X=\{X(t):-\tau\leq t\leq T\}$ which is a solution of \eqref{eq.maineqdiffnonexist}, \eqref{eq.icnonexist}.
\end{proposition}
We note that a solution does not exist \emph{for any} $T>0$.

It is the hypotheses $\epsilon=0$ that is crucial in ensuring the non--existence of a solution.
In \eqref{eq.maineqdiffnonexist} we may define the neutral functional $D$ by
\[
D(\varphi):=(1-\epsilon) \varphi(0)-\int_{-\tau}^0 w(s)h(\varphi(s))\,ds, \quad \varphi\in C([-\tau,0];\mathbb{R}).
\]
Suppose that $h$ is globally Lipschitz continuous with Lipschitz constant $k_h$. Let $\phi\in C([-\tau,0];\mathbb{R})$ and suppose that
$\varphi_1, \varphi_2 \in Q(\phi,s)$ for $s<\tau$. Clearly $D$ cannot be uniformly non--atomic at 0 on $C([-\tau,0];\mathbb{R})$ for otherwise
\eqref{eq.maineqdiffnonexist} would have a solution.

We now show, however, for $\epsilon\in(0,2)$ that $D$ is uniformly non--atomic at 0 on $C([-\tau,0];\mathbb{R})$, and so
\eqref{eq.maineqdiffnonexist} does have a solution. First note that
\[
D(\varphi_1)-D(\varphi_2)
=(1-\epsilon) (\varphi_1(0)-\varphi_2(0)) -\int_{-\tau}^0 w(u)\left(h(\varphi_1(u))-h(\varphi_2(u))\right)\,du.
\]
Since $\varphi_1(u)=\varphi_2(u)=\phi(u)$ for $u\in[-\tau,-s)$, we have
\begin{equation} \label{eq.Dnonatomicexamp1}
D(\varphi_1)-D(\varphi_2)
=(1-\epsilon) (\varphi_1(0)-\varphi_2(0)) % -\int_{-\tau}^{-s} w(u)\left(h(\varphi_1(u))-h(\varphi_2(u))\right)\,du -
-\int_{-s}^{0} w(u)\left(h(\varphi_1(u))-h(\varphi_2(u))\right)\,du.
\end{equation}
Therefore by \eqref{eq.Dnonatomicexamp1} we have
\begin{align*}
|D(\varphi_1)-D(\varphi_2)|
&\leq |1-\epsilon| |\varphi_1(0)-\varphi_2(0)| + \int_{-s}^{0} |w(u)||h(\varphi_1(u))-h(\varphi_2(u))|\,du\\
&\leq |1-\epsilon| |\varphi_1(0)-\varphi_2(0)| + k_h \int_{-s}^{0} |w(u)| |\varphi_1(u)-\varphi_2(u)|\,du\\
&\leq |1-\epsilon| \|\varphi_1-\varphi_2\|_{\text{sup}} + k_h \|\varphi_1-\varphi_2\|_{\text{sup}} \int_{-s}^{0} |w(u)| \,du\\
&%=\left( |1-\epsilon| + k_h \int_{-s}^{0} |w(u)| \,du \right)  \|\varphi_1-\varphi_2\|_{\text{sup}}
=\rho_0(s) \|\varphi_1-\varphi_2\|_{\text{sup}},
\end{align*}
where we define
\[
\rho_0(s):= |1-\epsilon| + k_h \int_{-s}^{0} |w(u)| \,du, \quad s\in[-\tau,0].
\]
Clearly $\rho_0$ is non--decreasing. For every $\epsilon\in (0,2)$ we have $|1-\epsilon|<1$,
so because $w$ is continuous, there exists a $T_1>0$ such that $\rho_0(T_1)<1$. In this case, $D$ is uniformly non--atomic at 0 on $C([-\tau,0];\mathbb{R})$. Therefore for $\epsilon\in(0,2)$ we see that \eqref{eq.maineqdiffnonexist} has a unique solution by Theorem \ref{thmexiuni}. In the case when $\epsilon>2$ or $\epsilon<0$, simply divide \eqref{eq.maineqdiffnonexist} by $\epsilon$. The properties on
$f$, $w$ and $h$ etc. guarantee the existence and uniqueness by Theorem \ref{thmexiuni} using the above arguments in the case $\epsilon=1$.
\begin{proposition} \label{prop.existstochadd}
Let $\tau>0$ and $\epsilon\neq 0$. Suppose $h\in C(\mathbb{R};\mathbb{R})$ is globally Lipschitz continuous, $w\in
C([-\tau,0];\mathbb{R})$, $\psi\in C([-\tau,0];\mathbb{R})$ and
$\sigma\neq 0$. Suppose also that there is $K>0$
\[
|f(\phi)-f(\varphi)|\leq K\sup_{-\tau\leq s\leq 0}|\phi(s)-\varphi(s)|, \quad \text{for all }\phi, \varphi\in C([-\tau,0];\mathbb{R})
\]
Let $T>0$. Then there is a unique solution $X=\{X(t):-\tau\leq t\leq T\}$ of \eqref{eq.maineqdiffnonexist}, \eqref{eq.icnonexist}.
\end{proposition}

\subsubsection{Equations with maximum functionals}
Let $\kappa>0$ and suppose that $g:C([-\tau,0];\mathbb{R})\to \mathbb{R}$ is globally Lipschitz continuous. Consider the SFDE
\begin{equation} \label{eq.sfdemaxtype}
d(X(t)+\kappa\max_{-\tau\leq s\leq 0} |X(t+s)|) = g(X_t)\,dB(t), \quad \text{$0\leq t\leq T$, a.s.}
\end{equation}
In the case when $\kappa\in(0,1)$, \eqref{eq.dlipcontract} holds for the functional $D$ defined by
\begin{equation} \label{eq.Dfunctmaxtypenonexist}
D(\varphi) = \kappa \max_{s\in [-\tau,0]} |\varphi(s)|, \quad \varphi \in C([-\tau,0];\mathbb{R}),
\end{equation}
and for any given $T>0$, \eqref{eq.sfdemaxtype} has a solution by Mao \cite[Theorem 6.2.2]{Mao97}. This could also be concluded from
the fact that $D$ is uniformly non--atomic at $0$ on $C([-\tau,0];\mathbb{R})$, in which case Theorem \ref{thmexiuni} applies.

We suppose now that $\kappa\geq 1$. We note that \eqref{eq.dlipcontract} does not apply to the functional $D$ in \eqref{eq.Dfunctmaxtypenonexist}. To see this consider
$\varphi_2\in C([-\tau,0],\mathbb{R})$ and let $\varphi_1=\alpha \varphi_2$ for some $\alpha>0$. Then
\begin{align*}
|D(\varphi_2)-D(\varphi_1)| &= |\kappa\|\varphi_2\|_{\text{sup}}- \kappa \|\varphi_1\|_{\text{sup}}| \\
&=\kappa|\|\varphi_2\|_{\text{sup}}- \alpha \|\varphi_2\|_{\text{sup}}|
=\kappa |1-\alpha| \|\varphi_2\|_{\text{sup}}.
\end{align*}
On the other hand $\kappa\|\varphi_2-\varphi_1\|_{\text{sup}}=\kappa\|\varphi_2 - \alpha\varphi_2\|_{\text{sup}}=\kappa|1- \alpha|\|\varphi_2 \|_{\text{sup}}$,
so
\[
|D(\varphi_2)-D(\varphi_1)|=\kappa\|\varphi_2-\varphi_1\|_{\text{sup}},
\]
which violates \eqref{eq.dlipcontract}, as $\kappa\geq 1$.

Also, we see that $D$ in \eqref{eq.Dfunctmaxtypenonexist} does not satisfy \eqref{eq.nonatomlip2a}. To see
this suppose that $\varphi_1,\varphi_2\in Q(s,0)$ is such that $\varphi_2(0)>0$, $\varphi_2$ is non--decreasing, and $\varphi_1=\alpha\varphi_2$
for $\alpha>0$. Then
\[
D(\varphi_2)=\kappa\max_{u\in[-\tau,0]} |\varphi_2(u)| = \kappa \max_{u\in[-s,0]} |\varphi_2(u)|=\kappa \max_{u\in[-s,0]} \varphi_2(u)=\kappa \varphi_2(0).
\]
Similarly
\[
D(\varphi_1)= \kappa \max_{u\in[-s,0]} |\varphi_1(u)|=\kappa \max_{u\in[-s,0]} \alpha\varphi_2(u)=\kappa\alpha \varphi_2(0).
\]
Hence $|D(\varphi_2)-D(\varphi_1)|=\kappa|1-\alpha|\varphi_2(0)$. On the other hand
\[
\|\varphi_2-\varphi_1\|_{\text{sup}} %= \max_{u\in[-\tau,0]}|\varphi_2(u)-\varphi_1(u)| =
=\max_{u\in[-s,0]}|\varphi_2(u)-\varphi_1(u)|= \max_{u\in[-s,0]} |1-\alpha||\varphi_2(u)|
=|1-\alpha|\varphi_2(0).
\]
Thus $|D(\varphi_2)-D(\varphi_1)|=\kappa \|\varphi_2-\varphi_1\|_{\text{sup}}$, so \eqref{eq.rho0rho0lt1} and \eqref{eq.nonatomlip2a} cannot both
be satisfied, because $\kappa\geq 1$.

We now prove that \eqref{eq.sfdemaxtype} does not have a solution.
\begin{proposition} \label{prop.sfdemaxtypenonexist}
Let $\tau>0$. Let $\psi\in C([-\tau,0];\mathbb{R})$ and
$\sigma\neq 0$. Suppose also that
\begin{equation} \label{eq.nondegdiffsfdemaxtype}
\text{There exists $\delta>0$ such that }
\delta:=\inf_{\varphi\in C([-\tau,0];\mathbb{R})} g^2(\varphi).
\end{equation}
Let $T>0$ and $\kappa\geq 1$. Then there is no process $X=\{X(t):-\tau\leq t\leq T\}$ which is a solution of \eqref{eq.sfdemaxtype}.
\end{proposition}
%%%%%%%%%%%%%%%%%%%%%%%%%%%%%%%%%%%%%%%%%%%%%%%%%%%%%%%%%%%%%%%%%%%%%%%%%%%%%%%%%%%%%%%%%%%%%%%%%%%%%%%%%%%%%%%%%%%%%%%%%%%%%%%%%%%%%%%%%%%%%
\section{Auxiliary Results}\label{sec:lemma}
The proofs of the main results are facilitated by a number of supporting lemmata. We state and discuss these here.

We first give a lemma which is necessary in proving the uniqueness and existence of the solution.

\begin{lemma}\label{lemafinitemo}
Let $X$ be the unique continuous solution of equation \eqref{eq.maineqdiff} with initial condition \eqref{eq.ic}. If both \eqref{eq.fglinboud} and \eqref{nint} hold, then for any $p\geq 2$, there exist positive constants $K_1$ and $K_2$ depending on $T$ such that
\begin{equation}\label{eq.lemafinitemo}
\mathbb{E}[\sup_{-\tau \leq s\leq t}|X(s)|^p]\leq K_1e^{K_2T}.
\end{equation}
\end{lemma}

In our proofs of moment estimates, we will need to use the fact that the $p$--th moment of the solution is a continuous function. Although the continuity of the moments is known for solutions of SNDEs, the contraction condition \eqref{eq.dlipcontract} is used in proving this continuity. Therefore, under our weaker assumptions, we need to prove this result afresh. To prove the continuity, we first need an elementary inequality.

\begin{lemma}  \label{lem.estcty}
Let $p\geq 1$. Suppose that $U,V\in \mathbb{R}^d$ are random variables in $L^{2(p-1)}$. If
$c_p>0$ is the number such that
\[
(a +b)^{2(p-1)} \leq c_p(a^{2(p-1)}+b^{2(p-1)}), \quad \text{for all $a,b\geq 0$},
\]
then
\[
\left|\mathbb{E}[|U|^p]-\mathbb{E}[|V|^p] \right|
\leq  p \left(c_p\mathbb{E}[(|U|^{2(p-1)}] +c_p\mathbb{E}[|V|^{2(p-1)}]\right)^{1/2}\mathbb{E}[|U-V|^2]^{1/2}.
\]
\end{lemma}

The continuity of the moments applies to general processes; since we will also employ it for an important auxiliary process, we do not confine
the scope of the result to the solution of \eqref{eq.maineqdiff}.

\begin{lemma} \label{lem.momcty}
Let $p\geq 1$. Let $\tau,T>0$. Let $X=\{X(t):t\in[-\tau,T]\}$ be a $\mathbb{R}^d$--valued stochastic process with a.s. continuous paths,
such that
\begin{equation} \label{eq.finite22pmin1mom}
\mathbb{E}[\max_{-\tau\leq s\leq T} |X(s)|^2]<+\infty, \quad
\mathbb{E}[\max_{-\tau\leq s\leq T} |X(s)|^{2(p-1)}]<+\infty.
\end{equation}
Then
\begin{equation} \label{eq.cty1}
\lim_{t\to s} \mathbb{E}[|X(t)-X(s)|^2]=0, \quad \text{for all $s\in[0,T]$},
\end{equation}
and so
\begin{equation} \label{eq.cty2}
\lim_{t\to s} \mathbb{E}[|X(t)|^p]=\mathbb{E}[|X(s)|^p] \quad \text{for all $s\in[0,T]$}.
\end{equation}
\end{lemma}

We find it useful to prove a variant of Gronwall's lemma. The
argument is a slight modification of arguments given in Gripenberg,
Londen and Staffans~\cite[Theorems 9.8.2 and 10.2.15]{GrLoSt90}. The result gives us
the freedom to \emph{construct} an upper bound via an integral
\emph{inequality}, rather than relying on precise knowledge of the
asymptotic behaviour of a solution of an \emph{equation}. We avail of this freedom
in proving a.s. and $p$-th mean exponential estimates on the solution of the neutral SFDE.
\begin{lemma}\label{lemagronwall}
Suppose that $\kappa\in M([0,\infty),\mathbb{R})$ is such that
$(-\kappa)$ has non--positive resolvent $\rho$ given by
\[
\rho+(-\kappa)\ast \rho = -\kappa.
\]
Let $f$ be in $L^1_{\text{loc}}(\mathbb{R}^+)$ and $x\in
L^1_{\text{loc}}(\mathbb{R}^+)$ obey
\begin{equation} \label{eq.ineqx}
x(t)\leq (\kappa\ast x)(t)+f(t), \quad t\geq 0.
\end{equation}
If $y\in L^1_{\text{loc}}(\mathbb{R}^+)$ obeys
\begin{equation} \label{eq.ineqy}
y(t)\geq (\kappa\ast y)(t)+f(t), \quad t\geq 0; \quad y(0)\geq x(0),
\end{equation}
then $x(t)\leq y(t)$ for all $t\geq 0$.
\end{lemma}

%%%%%%%%%%%%%%%%%%%%%%%%%%%%%%%%%%%%%%%%%%%%%%%%%%%%%%%%%%%%%%%%%%%%%%%%%%%%%%%%%%%%%%%%%%%%%%%%%%%%%%%%%%%%%%%%%%%%%%%%%%%%%%%%%%%%%%%%%%%%%%%%%%%
\section{Proof of Section \ref{sec:lemma}}
\subsection{Proof of Lemma \ref{lemafinitemo}}
First, consider $t\in [0, T_1]$. Define
\[
\xi_m:=T_1\wedge \inf\{t\in[0, T_1]\quad | \quad|X(t)|\geq m\}, \quad m\in \mathbb{N}.
\]
Set $X^m(t)=X(t\wedge \xi_m)$. Hence
\[
X^m(t)=\psi(0)-D(\psi)+D(X^m_t)+\int_0^tf(X^m_s)\,ds+\int_0^tg(X^m_s)\,dB(s).
\]
By the inequality (cf. \cite[Lemma 6.4.1]{Mao97}),
\begin{equation}\label{Maoineqp}
|a+b|^p\leq (1+\varepsilon^{\frac{1}{p-1}})^{p-1}(|a|^p+\frac{|b|^p}{\varepsilon}),
\quad \forall\,p>1,\,\, \varepsilon>0,\,\,\text{and} \,\,a,b\in\mathbb{R},
\end{equation}
it is easy to show that
\[
|X^m(t)|^p\leq (1+\varepsilon ^{\frac{1}{p-1}})^{p-1}\bigg(|D(X^m_t)-D(\psi)|^p+\frac{1}{\varepsilon}|J^m_1(t)|^p\bigg),
\]
where
\begin{equation}\label{varep}
0 < \varepsilon < \bigg(\frac{1}{k^{\frac{p}{3p-3}}}-1\bigg)^{p-1}
\end{equation}
$k$ is defined in \eqref{nint}, and
\[
J^m_1(t):=\psi(0)+\int_0^tf(X^m_s)\,ds+\int_0^tg(X^m_s)\,dB(s).
\]
Given \eqref{nint},and using \eqref{Maoineqp}, for any $\varepsilon>1$, we have
\begin{align*}
&|X^m(t)|^p
\\&\leq \frac{(1+\varepsilon^{\frac{1}{p-1}})^{2p-2}}{\varepsilon}|D(\psi)|^p+(1+\varepsilon^{\frac{1}{p-1}})^{2p-2}
|D(X_t^m)|^p+\frac{(1+\varepsilon^{\frac{1}{p-1}})^{p-1}}{\varepsilon}|J^m_1(t)|^p
\\&\leq \frac{(1+\varepsilon^{\frac{1}{p-1}})^{2p-2}}{\varepsilon}|D(\psi)|^p+\frac{(1+\varepsilon^{\frac{1}{p-1}})^{p-1}}{\varepsilon}|J^m_1(t)|^p
\\&\qquad+(1+\varepsilon^{\frac{1}{p-1}})^{2p-2}\bigg[K_D\bigg(1+\sup_{-\tau\leq s\leq -T_1}|X^m(t+s)|\bigg)+k\sup_{-T_1\leq s\leq 0}|X^m(t+s)|\bigg]^p
\\&\leq \frac{(1+\varepsilon^{\frac{1}{p-1}})^{2p-2}}{\varepsilon}|D(\psi)|^p+\frac{(1+\varepsilon^{\frac{1}{p-1}})^{p-1}}{\varepsilon}|J^m_1(t)|^p
\\&\qquad+(1+\varepsilon^{\frac{1}{p-1}})^{2p-2}[K_D+(K_D+k)\sup_{-\tau\leq s\leq 0}|\psi(s)|+k\sup_{0\leq s\leq t}|X^m(s)|]^p
\\&\leq \frac{(1+\varepsilon^{\frac{1}{p-1}})^{2p-2}}{\varepsilon}|D(\psi)|^p+\frac{(1+\varepsilon^{\frac{1}{p-1}})^{p-1}}{\varepsilon}|J^m_1(t)|^p
\\&\qquad+(1+\varepsilon^{\frac{1}{p-1}})^{3p-3}k^p\sup_{0\leq s\leq t}|X^m(s)|^p
\\&\qquad+\frac{(1+\varepsilon^{\frac{1}{p-1}})^{3p-3}}{\varepsilon}[K_D+(K_D+k)\sup_{-\tau\leq s\leq 0}|\psi(s)|]^p
\end{align*}
Thus
\begin{align*}
&\sup_{0\leq s\leq t}|X^m(s)|^p
\\&\leq \frac{(1+\varepsilon^{\frac{1}{p-1}})^{2p-2}}{\varepsilon}|D(\psi)|^p+
\frac{(1+\varepsilon^{\frac{1}{p-1}})^{3p-3}}{\varepsilon}[K_D+(K_D+k)\sup_{-\tau\leq s\leq 0}|\psi(s)|]^p
\\&+(1+\varepsilon^{\frac{1}{p-1}})^{3p-3}k^p\sup_{0\leq s\leq t}|X^m(s)|^p+\frac{(1+\varepsilon^{\frac{1}{p-1}})^{p-1}}{\varepsilon}\sup_{0\leq s\leq t}|J^m_1(t)|^p.
\end{align*}
Due to \eqref{varep}, $(1+\varepsilon^{\frac{1}{p-1}})^{3p-3}k^p<1$, the above inequality implies
\begin{multline*}
\sup_{0\leq s\leq t}|X^m(s)|^p\leq \frac{1}{1-(1+\varepsilon^{\frac{1}{p-1}})^{3p-3}k^p}\bigg\{\frac{(1+\varepsilon^{\frac{1}{p-1}})^{2p-2}}{\varepsilon}|D(\psi)|^p\\+
\frac{(1+\varepsilon^{\frac{1}{p-1}})^{3p-3}}{\varepsilon}[K_D+(K_D+k)\sup_{-\tau\leq s\leq 0}|\psi(s)|]^p\bigg\}
\\+\frac{(1+\varepsilon^{\frac{1}{p-1}})^{p-1}}{\varepsilon[1-(1+\varepsilon^{\frac{1}{p-1}})^{3p-3}k^p]}\sup_{0\leq s\leq t}|J^m_1(t)|^p.
\end{multline*}
Since
\[
\sup_{-\tau\leq s\leq t}|X^m(s)|^p\leq \sup_{-\tau\leq s\leq 0}|\psi(s)|^p+\sup_{0\leq s\leq t}|X^m(s)|^p,
\]
we get
\begin{multline}\label{supX}
\sup_{-\tau\leq s\leq t}|X^m(s)|^p\leq \bigg\{\frac{1}{1-(1+\varepsilon^{\frac{1}{p-1}})^{3p-3}k^p}\bigg[\frac{(1+\varepsilon^{\frac{1}{p-1}})^{2p-2}}{\varepsilon}|D(\psi)|^p\\+
\frac{(1+\varepsilon^{\frac{1}{p-1}})^{3p-3}}{\varepsilon}[K_D+(K_D+k)\sup_{-\tau\leq s\leq 0}|\psi(s)|]^p\bigg]
 +\sup_{-\tau\leq s\leq 0}|\psi(s)|^p\bigg\}\\+\frac{(1+\varepsilon^{\frac{1}{p-1}})^{p-1}}{\varepsilon[1-(1+\varepsilon^{\frac{1}{p-1}})^{3p-3}k^p]}\sup_{0\leq s\leq t}|J^m_1(t)|^p.
\end{multline}
Now
\begin{align*}
&\sup_{0\leq s\leq t}|J^m_1(s)|^p \\&= \sup_{0\leq s\leq t}\bigg|\psi(0)+\int_0^sf(X^m_u)\,du+\int_0^sg(X^m_u)\,dB(u)\bigg|^p
\\&\leq \frac{(1+\varepsilon^{\frac{1}{p-1}})^{p-1}}{\varepsilon}\sup_{-\tau\leq s\leq 0}|\psi(s)|^p
\\& \quad\quad\quad\quad
+(1+\varepsilon^{\frac{1}{p-1}})^{p-1}\sup_{0\leq s\leq t}\bigg|\int_0^sf(X^m_u)\,du+\int_0^sg(X^m_u)\,dB(u)\bigg|^p
\\&\leq \frac{(1+\varepsilon^{\frac{1}{p-1}})^{p-1}}{\varepsilon}\sup_{-\tau\leq s\leq 0}|\psi(s)|^p
+\frac{(1+\varepsilon^{\frac{1}{p-1}})^{2p-2}}{\varepsilon}\sup_{0\leq s\leq t}\bigg(\int_0^s\bigg|f(X^m_u)\bigg|\,du\bigg)^p+
\\&\quad\quad\quad\quad
+(1+\varepsilon^{\frac{1}{p-1}})^{2p-2}\sup_{0\leq s\leq t}\bigg|\int_0^sg(X^m_u)\,dB(u)\bigg|^p
\end{align*}
Taking expectations on both sides of the inequality, and let $\alpha=\varepsilon^{1/(p-1)}$, by Assumption \ref{ass.fgliplinbdd}, we have
\begin{align*}
\mathbb{E}[\sup_{0\leq s\leq t}|J^m_1(s)|^p]&\leq \bigg(\frac{1+\alpha}{\alpha}\bigg)^{p-1}\sup_{-\tau\leq s\leq 0}|\psi(s)|^p
\\&\quad\quad\quad+\frac{(1+\alpha)^{2p-2}}{\alpha^{p-1}}\mathbb{E}\bigg[\sup_{0\leq s\leq t}\bigg(\int_0^s\bar{K}(1+||X^m_u||_\text{sup})\,du\bigg)^p\bigg]
\\&\quad\quad\quad+(1+\alpha)^{2p-2}\mathbb{E}\bigg[\sup_{0\leq s\leq t}\bigg|\int_0^sg(X^m_u)\,dB(u)\bigg|^p\bigg].
\end{align*}
By the Burkholder-Davis-Gundy inequality, let $C_p:=[p^{p+1}/(2(p-1)^{p-1})]^{p/2}$, the above inequality implies that
\begin{align*}
\mathbb{E}[\sup_{0\leq s\leq t}|J^m_1(s)|^p]&\leq \bigg(\frac{1+\alpha}{\alpha}\bigg)^{p-1}\sup_{-\tau\leq s\leq 0}|\psi(s)|^p
\\&\quad\quad\quad+\frac{(1+\alpha)^{2p-2}}{\alpha^{p-1}}\bar{K}^p\mathbb{E}\bigg[\bigg(\int_0^t(1+||X^m_u||_\text{sup})\,du\bigg)^p\bigg]
\\&\quad\quad\quad+(1+\alpha)^{2p-2}C_p\mathbb{E}\bigg[\bigg(\int_0^t||g(X^m_s)||^2_s\,ds\bigg)^{\frac{p}{2}}\bigg]
\\&\leq \bigg(\frac{1+\alpha}{\alpha}\bigg)^{p-1}\sup_{-\tau\leq s\leq 0}|\psi(s)|^p
\\&\quad\quad\quad+\frac{(1+\alpha)^{2p-2}}{\alpha^{p-1}}\bar{K}^pT_1^{p-1}\mathbb{E}\bigg[\int_0^t(1+||X^m_u||_\text{sup})^p\,du\bigg]
\\&\quad\quad\quad+(1+\alpha)^{2p-2}C_p\bar{K}^p\mathbb{E}\bigg[\bigg(\int_0^t(1+||X^m_u||_\text{sup})^2\,du\bigg)^{\frac{p}{2}}\bigg]
\end{align*}
where we have used H\"{o}lder's inequality in the second line. Thus
\[
(1+||X^m_u||_\text{sup})^p\leq (1+\alpha)^{p-1}(\alpha^{1-p}+||X^m_u||_\text{sup}^p),
\]
and
\begin{align*}
\bigg(\int_0^t(1+||X^m_u||_\text{sup})^2\,du\bigg)^{\frac{p}{2}}
&\leq T_1^{\frac{(p-2)p}{4}}\int_0^t(1+||X^m_u||_\text{sup})^p\,du
\\& \leq (1+\alpha)^{p-1}T_1^{\frac{(p-2)p}{4}}\int_0^t(\alpha^{1-p}+||X^m_u||_\text{sup}^p)\,du.
\end{align*}
Hence
\begin{multline}
\mathbb{E}[\sup_{0\leq s\leq t}|J^m_1(s)|^p]\leq \bigg(\frac{1+\alpha}{\alpha}\bigg)^{p-1}\sup_{-\tau\leq s\leq 0}|\psi(s)|^p
\\+\bigg[\frac{(1+\alpha)^{3p-3}}{\alpha^{p-1}}\bar{K}^pT_1^{p-1}+(1+\alpha)^{3p-3}C_p\bar{K}^pT_1^{\frac{(p-2)p}{4}}\bigg]
\mathbb{E}\bigg[\int_0^t(\alpha^{1-p}+||X^m_u||^p_\text{sup})\,du\bigg].
\end{multline}
Taking expectations on both sides of \eqref{supX}, and inserting the above inequality into \eqref{supX}, we have
\begin{align*}
\frac{1}{\varepsilon}+\mathbb{E}[\sup_{-\tau\leq s\leq t}|X^m(s)|^p]&\leq \kappa_1+\kappa_2\int_0^t\bigg(\frac{1}{\varepsilon}+\mathbb{E}[||X^m_u||^p_\text{sup}]\bigg)\,du
\\&\leq (\frac{1}{\varepsilon}+\kappa_1)+\kappa_2\int_0^t\bigg(\frac{1}{\varepsilon}+\mathbb{E}[\sup_{-\tau\leq u\leq s}|X^m(u)|^p\bigg)\,du,
\end{align*}
where
\begin{multline*}
\kappa_1:=\frac{1}{1-(1+\varepsilon^{\frac{1}{p-1}})^{3p-3}k^p}\bigg[\frac{(1+\varepsilon^{\frac{1}{p-1}})^{2p-2}}{\varepsilon}|D(\psi)|^p\\+
\frac{(1+\varepsilon^{\frac{1}{p-1}})^{3p-3}}{\varepsilon}[K_D+(K_D+k)\sup_{-\tau\leq s\leq 0}|\psi(s)|]^p\bigg]
 \\+\bigg[1+\frac{(1+\varepsilon^{\frac{1}{p-1}})^{p-1}}{\varepsilon}\bigg]\sup_{-\tau\leq s\leq 0}|\psi(s)|^p,
\end{multline*}
and
\begin{multline*}
\kappa_2:=\frac{(1+\varepsilon^{\frac{1}{p-1}})^{p-1}}{\varepsilon[1-(1+\varepsilon^{\frac{1}{p-1}})^{3p-3}k^p]}\times
\\ \bigg[\frac{(1+\varepsilon^{\frac{1}{p-1}})^{3p-3}}{\varepsilon}\bar{K}^pT_1^{p-1}+(1+\varepsilon^{\frac{1}{p-1}})
^{3p-3}C_p\bar{K}^pT_1^{\frac{(p-2)p}{4}}\bigg].
\end{multline*}
Now the Gronwall inequality yields that
\[
\frac{1}{\varepsilon}+\mathbb{E}[\sup_{-\tau\leq s\leq T_1}|X^m(s)|^p]\leq (\frac{1}{\varepsilon}+\kappa_1)e^{\kappa_2T_1},
\]
Consequently
\[
\mathbb{E}[\sup_{-\tau\leq s\leq T_1}|X^m(s)|^p]\leq (\frac{1}{\varepsilon}+\kappa_1)e^{\kappa_2T_1}.
\]
Letting $m\to\infty$ and $\varepsilon\to [1/k^{p/(3p-3)}-1]^{p-1}$, we get
\[
\mathbb{E}[\sup_{-\tau\leq s\leq T_1}|X(s)|^p]\leq \bigg[\bigg(\frac{1}{k^{\frac{p}{3p-3}}}-1\bigg)^{p-1}+\kappa_1\bigg]e^{\kappa_2T_1}.
\]
For $t\in[nT_1, (n+1)T_1]$ ($n\in\mathbb{N}$), assertion \eqref{eq.lemafinitemo} can be shown by
applying the same analysis as in the case of $t\in [0, T_1]$.

\subsection{Proof of Lemma \ref{lem.estcty}}
Let $x,y\geq 0$ and $p\geq 1$. Then there exists $\theta(x,y)\in[0,1]$ such that
\[
x^p-y^p = p [\theta x + (1-\theta) y]^{p-1}(x-y).
\]
Thus for $U$, $V\in\mathbb{R}^d$ we have $\theta(U,V)\in [0,1]$ such that
\[
|U|^p-|V|^p = p [\theta |U| + (1-\theta) |V|]^{p-1}(|U|-|V|).
\]
Therefore
\begin{align*}
\mathbb{E}[|U|^p]-\mathbb{E}[|V|^p]
&= p \mathbb{E}[[\theta |U| + (1-\theta) |V|]^{p-1}(|U|-|V|)]\\
&\leq p \mathbb{E}[[\theta |U| + (1-\theta) |V|]^{2(p-1)}]^{1/2}\mathbb{E}[(|U|-|V|)^2]^{1/2}\\
&\leq p \mathbb{E}[[|U| + |V|]^{2(p-1)}]^{1/2}\mathbb{E}[(|U|-|V|)^2]^{1/2}.
\end{align*}
Similarly, as $|V|^p-|U|^p= p [\theta |U| + (1-\theta) |V|]^{p-1}(|V|-|U|)$, we have
\begin{align*}
\mathbb{E}[|V|^p]-\mathbb{E}[|U|^p]
&= p \mathbb{E}[[\theta |U| + (1-\theta) |V|]^{p-1}(|V|-|U|)]\\
&\leq p \mathbb{E}[[\theta |U| + (1-\theta) |V|]^{2(p-1)}]^{1/2}\mathbb{E}[(|V|-|U|)^2]^{1/2}\\
&=p \mathbb{E}[[\theta |U| + (1-\theta) |V|]^{2(p-1)}]^{1/2}\mathbb{E}[(|U|-|V|)^2]^{1/2}\\
&\leq p \mathbb{E}[[|U| +|V|]^{2(p-1)}]^{1/2}\mathbb{E}[(|U|-|V|)^2]^{1/2}.
\end{align*}
Therefore
\begin{align*}
\left|\mathbb{E}[|U|^p]-\mathbb{E}[|V|^p] \right|
&\leq p \mathbb{E}[[|U| +|V|]^{2(p-1)}]^{1/2}\mathbb{E}[(|U|-|V|)^2]^{1/2}\\
&= p \mathbb{E}[[|U| +|V|]^{2(p-1)}]^{1/2}\mathbb{E}[||U|-|V||^2]^{1/2}
\end{align*}
Now $||U|-|V||\leq |U-V|$, so $||U|-|V||^2\leq |U-V|^2$. Therefore
\[
\left|\mathbb{E}[|U|^p]-\mathbb{E}[|V|^p] \right|
\leq  p \mathbb{E}[[|U| +|V|]^{2(p-1)}]^{1/2}\mathbb{E}[|U-V|^2]^{1/2}.
\]
Since $(a+b)^{2(p-1)}\leq c_p(a^{2(p-1)}+b^{2(p-1)})$ for all $a,b\geq 0$, we have
\[
\left|\mathbb{E}[|U|^p]-\mathbb{E}[|V|^p] \right|
\leq  p \left(c_p\mathbb{E}[(|U|^{2(p-1)}] +c_p\mathbb{E}[|V|^{2(p-1)}]\right)^{1/2}\mathbb{E}[|U-V|^2]^{1/2},
\]
as required.

\subsection{Proof of Lemma~\ref{lem.momcty}}
Let $0\leq s\leq t\leq T$. We first prove \eqref{eq.cty1}. By the continuity of the sample paths, we have $\lim_{t\to s}X(t)=X(s)$ a.s. for
each $s\in[0,T]$. On the other hand, because
\[
|X(t)|\leq \max_{0\leq u\leq T} |X(u)|,
\]
we have that $|X(t)|$ is dominated by a random variable which is in $L^2$ by \eqref{eq.finite22pmin1mom}.
Then by the Dominated Convergence Theorem, we have that $X(t)$ converges to $X(s)$ in $L^2$ viz.,
\[
\lim_{t\to s}\mathbb{E}[|X(t)-X(s)|^2]=0,
\]
which is \eqref{eq.cty1}. Now we prove \eqref{eq.cty2}. Let $0\leq s\leq t\leq T$.
 Define $M_p(T):=\mathbb{E}[\max_{-\tau\leq s\leq T} |X(s)|^{2(p-1)}]$. Since \eqref{eq.finite22pmin1mom} holds, by Lemma~\ref{lem.estcty}
\begin{align*}
\lefteqn{\left|\mathbb{E}[|X(t)|^p]-\mathbb{E}[|X(s)|^p] \right|}\\
&\leq  p \left(c_p\mathbb{E}[(|X(t)|^{2(p-1)}] +c_p\mathbb{E}[|X(s)|^{2(p-1)}]\right)^{1/2}\mathbb{E}[|X(t)-X(s)|^2]^{1/2}\\
&\leq  p \left(2c_p M_p(T)\right)^{1/2}\mathbb{E}[|X(t)-X(s)|^2]^{1/2}.
\end{align*}
Now \eqref{eq.cty1} implies \eqref{eq.cty2}.

\subsection{Proof of Lemma \ref{lemagronwall}}
By \eqref{eq.ineqx} and \eqref{eq.ineqy}, there are $g\geq 0$ and a
$h\geq 0$, both in $L^1_{\text{loc}}(\mathbb{R}^+)$ such that
\[
x(t)= (\kappa\ast x)(t)+f(t)-g(t), \quad y(t)=(\kappa\ast
y)(t)+f(t)+h(t), \quad t\geq 0.
\]
Since $\rho$ is the resolvent of $-\kappa$, we have the variation of
constants formulae:
\[
x=f-g-\rho\ast(f-g), \quad y=f+h-\rho\ast(f+h).
\]
Therefore
\[
\kappa\ast x=\kappa \ast (f-g)-\kappa\ast\rho\ast(f-g)
=[\kappa-\kappa\ast\rho]\ast f- [\kappa-\kappa\ast\rho]\ast g=-\rho\ast f+\rho\ast g.
\]
Similarly $\kappa\ast y=-\rho\ast f-\rho\ast h$. Hence
\[
x(t) \leq (\kappa\ast x)(t)+f(t)= -(\rho\ast f)(t)+(\rho\ast g)(t)+f(t)\leq f(t) -(\rho\ast f)(t),
\]
where we have used the fact that $g$ is non--negative and $\rho$ is non--positive at the last step. Similarly
\[
y(t) \geq (\kappa\ast y)(t)+f(t)  -(\rho\ast f)(t)-(\rho\ast h)(t)+f(t)\geq f(t) -(\rho\ast f)(t),
\]
where we have used the fact that $h$ is non--negative and $\rho$ is non--positive at the last step. Therefore
$x(t) \leq  f(t) -(\rho\ast f)(t)\leq y(t)$ for all $t\geq 0$, which proves the claim.

%%%%%%%%%%%%%%%%%%%%%%%%%%%%%%%%%%%%%%%%%%%%%%%%%%%%%%%%%%%%%%%%%%%%%%%%%%%%%%%%%%%%%%%%%%%%%%%%%%%%%%%%%%%%%%%%%%%%%%%%%%%%%%%%%%%%%%%%%%%%%%%%%%%
\section{Proofs of Section \ref{sec:existence}}
\subsection{Proof of Theorem \ref{thmexiuni}}
We first establish the existence of the solution on $[0,T_1]$, where $T_1\in(0,\delta)$ as defined in Assumption \ref{assD}.
Define that for $n=0, 1, 2, ...,$ $X_{1,0}^n=\psi$ and $X_1^0(t)=\psi(0)$ for $0\leq t\leq T_1$. Define the Picard Iteration, for $n\in \mathbb{N}$, $t\in[0, T_1]$,
\begin{equation}\label{picard}
X_1^n(t)-D(X_{1,t}^{n-1})=\psi(0)-D(\psi)+\int_0^tf(X_{1,s}^{n-1})\,ds+\int_0^tg(X_{1,s}^{n-1})\,dB(s).
\end{equation} Hence
\[
X_1^1(t)-X_1^0(t)=D(X^0_{1,t})-D(\psi)+\int_0^tf(X_{1,s}^0)\,ds+\int_0^tg(X_{1,s}^0)\,dB(s).
\]
By Assumption \ref{assDlineargr},
\begin{align*}
|X_1^1(t)-X_1^0(t)|^2&\leq \frac{1}{\alpha}|D(X_{1,t}^0)-D(\psi)|^2+\frac{1}{1-\alpha}|I(t)|^2
\\ &\leq \frac{1}{\alpha}\bigg(K_D(1+\sup_{-\tau\leq s \leq -T_1}|X_1^0(t+s)|)\\&\quad\quad\quad\quad\quad
+k\sup_{-T_1\leq s\leq 0}|X_1^0(t+s)|+|D(\psi)|\bigg)^2
\end{align*}
where
\[
I(t):=\int_0^tf(X^0_{1,s})\,ds+\int_0^tg(X^0_{1,s})\,dB(s).
\]
It follows that
\begin{align*}
&\sup_{0\leq t \leq T_1}|X_1^1(t)-X_1^0(t)|^2
\\&\leq \frac{1}{\alpha}\bigg(K_D(1+\sup_{-\tau\leq s\leq 0}|\psi(s)|)+k\sup_{-T_1\leq s\leq T_1}|X_1^0(s)|
\\&\quad\quad\quad\quad\quad\quad\quad\quad\quad\quad\quad\quad\quad+|D(\psi)|\bigg)^2+\frac{1}{1-\alpha}\sup_{0\leq s\leq T_1}|I(t)|^2
\\&=\frac{1}{\alpha}\bigg(K_D+(K_D+k)\sup_{-\tau\leq s\leq 0}|\psi(s)|
+|D(\psi)|\bigg)^2+\frac{1}{1-\alpha}\sup_{0\leq s\leq T_1}|I(t)|^2
\end{align*}
By Assumption \ref{ass.fgliplinbdd}, it can be shown that
\[
\mathbb{E}\bigg[\sup_{0\leq t\leq T_1}|I(t)|^2\bigg]\leq 2\bar{K}T_1(T_1+4)(\sup_{-T_1\leq s\leq 0}|\psi(s)|^2+1).
\]
This implies that
\begin{multline}\label{expC}
\mathbb{E}\bigg[\sup_{0\leq t\leq T_1}|X_1^1(t)-X_1^0(t)|^2\bigg]\leq
\frac{1}{\alpha}\bigg(K_D+(K_D+k)\sup_{-\tau\leq s\leq 0}|\psi(s)|
+|D(\psi)|\bigg)^2
\\+
\frac{2\bar{K}T_1(T_1+4)}{1-\alpha}(\sup_{-T_1\leq s\leq 0}|\psi(s)|^2+1|)=:C.
\end{multline}
Now for all $n\in\mathbb{N}$ and $0\leq t\leq T_1<\delta$ ($\delta$ is defined in Assumption \ref{assD}), follow the same argument as in the proof of the uniqueness, we have $D_0(X^n_{1,t})-D_0(X^{n-1}_{1,t})=0$.
Therefore
\begin{multline*}
X_1^{n+1}(t)-X_1^n(t)=D_1(X_{1,t}^n)-D_1(X_{1,t}^{n-1})\\
+\int_0^t\big(f(X^n_{1,s})-f(X^{n-1}_{1,s})\big)\,ds+\int_0^t\big(g(X^n_{1,s})-g(X^{n-1}_{1,s})\big)\,dB(s).
\end{multline*}
Again by \eqref{eq.nonatomlip2}, we have
\begin{align*}
&|D_1(X^n_{1,t})-D_1(X^{n-1}_{1,t})|
\\&\leq k\|X_{1,t}^n-X_{1,t}^{n-1}\|_{\text{sup}}
\\&= k\max\{\sup_{-\tau\leq s\leq -T_1}|X_1^n(t+s)-X_1^{n-1}(t+s)|,
\\&\qquad\qquad\qquad\qquad\qquad\qquad\qquad\sup_{-T_1\leq s\leq 0}|X_1^n(t+s)-X_1^{n-1}(t+s)|\}
\\&= k\sup_{-T_1\leq s\leq 0}|X_1^n(t+s)-X_1^{n-1}(t+s)|\\
&=k \sup_{0\leq s\leq t}|X_1^n(s)-X_1^{n-1}(s)|.
\end{align*}
Apply the same analysis as in the proof of the uniqueness, we get
\begin{align}\label{expgamma}
&\mathbb{E}\bigg[\sup_{0\leq t\leq T_1}|X_1^{n+1}(t)-X_1^n(t)|^2\bigg]
\\ \nonumber &\leq \frac{k^2}{\alpha}\mathbb{E}\bigg[\sup_{0\leq t\leq T_1}|X_1^n(t)-X_1^{n-1}(t)|^2\bigg]
\\ \nonumber &\quad\quad\quad\quad\quad\quad\quad\quad\quad+\frac{2K(T_1+4)}{1-\alpha}\int_0^{T_1}\mathbb{E}\bigg[\sup_{0\leq s\leq t}|X_1^n(s)-X_1^{n-1}(s)|^2\bigg]\,dt
\\ \nonumber &\leq \bigg(\frac{k^2}{\alpha}+\frac{2KT_1(T_1+4)}{1-\alpha}\bigg)\mathbb{E}\bigg[\sup_{0\leq t\leq T_1}|X_1^n(s)-X_1^{n-1}(s)|^2\bigg].
\end{align}
Now let
\[
\gamma:= \frac{k^2}{\alpha}+\frac{2KT_1(T_1+4)}{1-\alpha}.
\]
We show that there exist such $T_1$ and $\alpha$ so that $\gamma<1$.
Fix $0<\mu<1$. Choose $T_1$ such that $k=\rho_0(T_1)<\mu$ and $2KT_1(T_1+4)<(1-\mu^2)^2/[2(1+\mu^2)]$. Let $\alpha=(1/2)\mu^2+(1/2)$, then
$k^2<\mu^2<\alpha<1$, which implies $\gamma<1$.
Combining \eqref{expgamma} with \eqref{expC}, we have
\begin{equation}\label{L2converge} \mathbb{E}\bigg[\sup_{0\leq t\leq T_1}|X_1^{n+1}(t)-X_1^n(t)|^2\bigg]\leq \gamma^nC.
\end{equation} Choose $\epsilon>0$, so that $(1+\epsilon)\gamma<1$. Hence by Chebyshev's inequality,
\[
\mathbb{P}\bigg\{\sup_{0\leq t\leq T_1}|X_1^{n+1}(t)-X_1^n(t)|>\frac{1}{(1+\epsilon)^n}\bigg\}\leq (1+\epsilon)^n\gamma^nC.
\]
Since $\sum_{n=0}^\infty(1+\epsilon)^n\gamma^nC<\infty$, by Borel-Cantelli lemma, for almost all $\omega\in\Omega$, there exists $n_0=n_0(\omega)\in\mathbb{N}$
such that
\[
\sup_{0\leq t\leq T_1}|X_1^{n+1}(t)-X_1^n(t)|\leq \frac{1}{(1+\epsilon)^n}, \quad\text{for  } n>n_0.
\]
This implies that
\[
X_1^n(t)=X_1^0(t)+\sum_{i=0}^{n-1}[X_1^{n-1}(t)-X_1^n(t)],
\]
converge uniformly on $t\in [0, T_1]$ a.s. Let the limit be $X_1(t)$ for $t\in [0, T_1]$ which is continuous and $\mathcal{F}(t)$-adapted.
Moreover, by \eqref{L2converge}, $\{X_1^n(t)\}_{n\in \mathbb{N}}\to X_1(t)$ in $L^2$ on $t\in [0, T_1]$. By Lemma \ref{lemafinitemo}, $X_1(\cdot)\in \mathcal{M}^2([-\tau, T_1];\mathbb{R}^d)$. Note that
\begin{align*}
\mathbb{E}\bigg[\bigg|\int_0^tf(X_{1,s}^n)\,ds-\int_0^tf(X_{1,s})\,ds\bigg|^2\bigg]
& \leq \mathbb{E}\bigg[\bigg(\int_0^t|f(X_{1,s}^n)-f(X_{1,s})|\,ds\bigg)^2\bigg]
\\& \leq \mathbb{E}\bigg[\bigg(\int_0^tK\|X_{1,s}^n-X_{1,s}\|_{\text{sup}}\,ds\bigg)^2\bigg]
\\& \leq K^2 T_1^2 \int_0^{T_1}\mathbb{E}[\|X_{1,s}^n-X_{1,s}\|_{\text{sup}}^2]\,ds,
\\&\to 0,\quad\text{as }\,n\to\infty,
\end{align*}
and
\begin{align*}
\mathbb{E}\bigg[\bigg|\int_0^tg(X_{1,s}^n)\,dB(s)&-\int_0^tg(X_{1,s})\,dB(s)\bigg|^2\bigg]
\\&= \mathbb{E}\bigg[\bigg|\int_0^t\bigg(g(X_{1,s}^n)-g(X_{1,s})\bigg)\,dB(s)\bigg|^2\bigg]
\\&= \mathbb{E}\bigg[\int_0^t\big|g(X_{1,s}^n)-g(X_{1,s})\big|^2\,ds\bigg]
\\&\leq K^2\int_0^{T_1}\mathbb{E}[\|X_{1,s}^n-X_{1,s}\|_{\text{sup}}^2]\,ds
\\&\to 0,\quad\text{as }\,n\to\infty,
\end{align*}
and
\[
\mathbb{E}[|D(X_{1,t}^n)-D(X_{1,t})|]\leq k\mathbb{E}[\|X_{1,t}^n-X_{1,t}\|]\to 0,\quad \text{as } n \to \infty.
\]
Hence let $n\to\infty$ in \eqref{picard}, almost surely that
\[
X_1(t)=\psi(0)-D(\psi)+D(X_{1,t})+\int_0^tf(X_{1,s})\,ds+\int_0^tg(X_{1,s})\,dB(s).
\]
Therefore $\{X_1(t)\}_{t\in[0, T_1]}$ is the solution on $[0,T_1]$ on an almost sure event $\Omega_{T_1}$.
We now prove the existence of the solution on the interval $[T_1, 2T_1]$. Define $X_{2,T_1}^n=X_{1,T_1}$ for $n=0,1,2...$, and $X_2^0(t)=X_1(T_1)$ for $t\in [T_1, 2T_1]$. Define the Picard Iteration, for $ n\in\mathbb{N}$,
\[
X_2^n(t)-D(X_{2,t}^{n-1})=X_1(T_1)-D(X_{1,T_1})+\int_{T_1}^tf(X_{2,s}^{n-1})\,ds+\int_{T_1}^tg(X_{2,s}^{n-1})\,dB(s).
\]
Following the same argument as in the case of $t\in [0, T_1]$, it can be shown that there exists continuous $\{X_2(t)\}_{t\in[T_1,2T_1]}$ such that $X_2^n(t)\to X_2(t)$ in $L^2$ for $t\in [T_1,2T_1]$ almost surely. Moreover, $X_2(\cdot)\in \mathcal{M}^2([T_1, 2T_1];\mathbb{R}^d)$, and $X_2(\cdot)$ almost surely satisfies the equation
\[
X_2(t)=X_1(T_1)-D(X_{1,T_1})+D(X_{2,t})+\int_{T_1}^tf(X_{2,s})\,ds+\int_{T_1}^tg(X_{2,s})\,dB(s).
\]
Therefore $\{X_2(t)\}_{t\in [T_1, 2T_1]}$ is the solution on $[T_1, 2T_1]$ on an almost sure event $\Omega_{2T_1}$. Let $X(t):=\{X_n(t)\cdot I_{\{t\in [nT_1, (n+1)T_1]\}}\}_{n\in \mathbb{N}\cup \{0\}}$, then $X(\cdot)$ is the solution of \eqref{eq:main} on the entire interval $[0, T]$ which is in $\mathcal{M}^2([0, T];\mathbb{R})$.

For the \emph{uniqueness}, consider $t\in [0, T_1]$, suppose that both $X$ and $Y$ are solutions to \eqref{eq:main}, with initial solution
$X(t)=Y(t)=\psi(t)$ for $t\in [-\tau,0]$. Then
\begin{multline*}
X(t)-Y(t)=D_0(X_t)-D_0(Y_t)+D_1(X_t)-D_1(Y_t)+\int_0^t \big(f(X_s)-f(Y_s)\big)\,ds\\+\int_0^t \big(g(X_s)-g(Y_s)\big)\,dB(s).
\end{multline*}
Let $s\in[-\tau,-\delta]$, by \eqref{D0}, we have $t+s\leq T_1-\delta<0$, and so
$X(t+s)=Y(t+s)=\psi(t+s)$. Then $|D_0(X_t)-D_0(Y_t)|=0$. Hence
\begin{multline*}
|X(t)-Y(t)|\leq |D_1(X_t)-D_1(Y_t)|\\+\bigg|\int_0^t \big(f(X_s)-f(Y_s)\big)\,ds+\int_0^t\big(g(X_s)-g(Y_s)\big)\,dB(s)\bigg|.
\end{multline*}
Let $k^2<\alpha<1$, where $k$ is given by \eqref{eq.rho0rho0lt1}. Then we get
\begin{equation*}
|X(t)-Y(t)|^2\leq \frac{1}{\alpha}|D_1(X_t)-D_1(Y_t)|^2+\frac{1}{1-\alpha}|J(t)|^2,
\end{equation*}
where we have used the inequality (cf. \cite[Lemma 6.2.3]{Mao97})
\begin{equation}\label{Maoineq}
(a+b)^2\leq \frac{1}{\alpha}a^2+\frac{1}{1-\alpha}b^2,\quad 0<\alpha<1.
\end{equation}
and define
\[
J(t):= \int_0^t \big(f(X_s)-f(Y_s)\big)\,ds+\int_0^t\big(g(X_s)-g(Y_s)\big)\,dB(s).
\]
Now by \eqref{eq.nonatomlip2}, since $0\leq t\leq T_1$,
\begin{align*}
&|D_1(X_t)-D_1(Y_t)|
\\&\leq k\|X_t-Y_t\|_{\text{sup}}
\\&= k\{\sup_{-\tau\leq s\leq -T_1}|X(t+s)-Y(t+s)|, \sup_{-T_1\leq s\leq 0}|X(t+s)-Y(t+s)|\}
\\&= k\sup_{-T_1\leq s\leq 0}|X(t+s)-Y(t+s)|.
\end{align*}
Therefore
\begin{align*}
|X(t)-Y(t)|^2&\leq \frac{k^2}{\alpha}\sup_{-T_1\leq s \leq 0}|X(t+s)-Y(t+s)|^2+\frac{1}{1-\alpha}|J(t)|^2
\\ &=\frac{k^2}{\alpha}\sup_{0\leq s\leq t}|X(s)-Y(s)|^2+\frac{1}{1-\alpha}|J(t)|^2.
\end{align*}
Moreover,
\[
\sup_{0\leq s \leq t}|X(s)-Y(s)|^2\leq \frac{k^2}{\alpha}\sup_{0\leq s\leq t}|X(s)-Y(s)|^2+\frac{1}{1-\alpha}\sup_{0\leq s\leq t}|J(t)|^2.
\]
Since $\alpha$ has been chosen such that $0<k^2<\alpha<1$, it follows that
\[
\sup_{0\leq s\leq t}|X(s)-Y(s)|^2\leq \frac{1}{(1-\alpha)(1-\frac{k^2}{\alpha})}\sup_{0\leq s\leq t}|J(t)|^2.
\]
Now, by \eqref{ass.fgliplinbdd} and similar argument as in the proof of Lemma \ref{lemafinitemo}, it is easy to show that
\[
\mathbb{E}\bigg[\sup_{0\leq s\leq t}|J(t)|^2\bigg]\leq 2K(T_1+4)\int_0^t\mathbb{E}\bigg[\sup_{0\leq u\leq s}|X(u)-Y(u)|^2\bigg]\,ds.
\]
It follows that
\[
\mathbb{E}\bigg[\sup_{0\leq s\leq t}|X(s)-Y(s)|^2\bigg]\leq \frac{2K(T_1+4)}{(1-\alpha)(1-\frac{k^2}{\alpha})}\int_0^t\mathbb{E}\bigg[\sup_{0\leq
u\leq s}|X(u)-Y(u)|^2\bigg]\,ds.
\]
Using Gronwall's inequality, we have that
\[
\forall\, 0\leq t\leq T_1, \quad \mathbb{E}\bigg[\sup_{0\leq s\leq t}|X(s)-Y(s)|^2\bigg]=0,
\]
which implies that
\[
\mathbb{E}\bigg[\sup_{0\leq t\leq T_1}|X(t)-Y(t)|^2\bigg]=0.
\]
Therefore we can conclude that on an a.s. event $\Omega_{T_1}$, for all $0\leq t\leq T_1$, $X(t)=Y(t)$ a.s. Apply the same argument on the interval
$[T_1, 2T_1]$ given $X(t)=Y(t)$ on $[-\tau, T_1]$ a.s., it can be shown that $X(t)=Y(t)$ on the entire interval $[-\tau, T]$ a.s.

\subsection{Proof of Theorem \ref{thmexpest}}
Let $Y(t):=X(t)-D(X_t)$, then by the inequality \eqref{Maoineqp},
we have
\begin{equation}\label{decompXp}
|X(t)|^p\leq (1+\varepsilon^{\frac{1}{p-1}})^{p-1}(|Y(t)|^p+\frac{1}{\varepsilon}|D(X_t)|^p).
\end{equation}
By It\^{o}'s formula,
\begin{multline*}
|Y(t)|^p=|\psi(0)-D(\psi)|^p+\int_0^t\bigg(p|Y(s)|^{p-2}Y^T(s)f(X_s)\\+\frac{p(p-1)}{2}|Y(s)|^{p-2}||g(X_s)||^2\,\bigg)ds
+\int_0^tp|Y(s)|^{p-2}Y^T(s)g(X_s)\,dB(s).
\end{multline*}
Hence if
\begin{equation}\label{finite4thmome}
\mathbb{E}\bigg[\int_0^t|Y(s)|^{2p-2}||g(X_s)||^2\,ds\bigg]<\infty,
\end{equation}
we get
\begin{multline*}
\mathbb{E}[|Y(t)|^p]=|\psi(0)-D(\psi)|^p+\mathbb{E}\bigg[\int_0^t\bigg(p|Y(s)|^{p-2}Y^T(s)f(X_s)\\
+\frac{p(p-1)}{2}|Y(s)|^{p-2}||g(X_s)||^2\,\bigg)ds\bigg].
\end{multline*}
We assume \eqref{finite4thmome} holds at the moment, and will show that it is true at the end of this proof.
Define $x(t):=\mathbb{E}[|X(t)|^p]$, and $y(t):=\mathbb{E}[|Y(t)|^p]$. Then
\begin{align*}
y(t+h)-y(t)&=\int_t^{t+h}\mathbb{E}\big[p|Y(s)|^{p-2}Y^T(s)f(X_s)+\frac{p(p-1)}{2}|Y(s)|^{p-2}||g(X_s)||^2\big]\,ds
\\&\leq \int_t^{t+h}\mathbb{E}\big[p|Y(s)|^{p-1}|f(X_s)|+\frac{p(p-1)}{2}|Y(s)|^{p-2}||g(X_s)||^2\big]\,ds
\\&\leq \int_t^{t+h}\bigg\{p\mathbb{E}\bigg[\frac{\varepsilon (p-1)}{p}|Y(s)|^p+\frac{|f(X_s)|^p}{p\varepsilon^{p-1}}\bigg]
\\&\quad\quad\quad\quad\quad\quad\quad\quad +\frac{p(p-1)}{2}\mathbb{E}\bigg[\frac{\varepsilon(p-2)}{p}|Y(s)|^p+\frac{2||g(X_s)||^p}{p\varepsilon^{(p-2)/2}}\bigg]\bigg\}\,ds
\\&=\int_t^{t+h}\bigg\{\frac{\varepsilon p(p-1)}{2} y(s)+\frac{1}{\varepsilon ^{p-1}}\mathbb{E}[|f(X_s)|^p]+\frac{p-1}{\varepsilon^{(p-2)/2}}
\mathbb{E}[||g(X_s)||^p]\bigg\}\,ds
\\&\leq \int_t^{t+h}\bigg\{\frac{\varepsilon p(p-1)}{2} y(s)+\frac{1}{\varepsilon ^{p-1}}\mathbb{E}\bigg[C_f+\int_{[-\tau,0]}\nu(du)|X(u+s)|^p\bigg]
\\&\quad\quad\quad\quad\quad\quad\quad\quad+\frac{p-1}{\varepsilon^{(p-2)/2}}\mathbb{E}\bigg[C_g+\int_{[-\tau,0]}\eta (du)|X(u+s)|^p\bigg]\bigg\}\,ds,
\end{align*}
where we have used the inequalities (cf. \cite[Lemma 6.2.4]{Mao97})
\[
\forall\, p\geq 2, \,\,\text{and}\,\,\varepsilon, a, b>0, \quad a^{p-1}b\leq \frac{\varepsilon(p-1)a^p}{p}+\frac{b^p}{p\varepsilon^{p-1}}
\]
and
\[
a^{p-2}b^2\leq \frac{\varepsilon(p-2)a^p}{p}+\frac{2b^p}{p\varepsilon^{(p-2)/2}},
\]
in the second inequality, conditions \eqref{f} and \eqref{g} in the last inequality. By the continuity of $t \mapsto \mathbb{E}[|X(t)|^p]$ and
$t \mapsto \mathbb{E}[|Y(t)|^p]$, it is then easy to see that
\[
D^+ y(t)\leq \frac{\varepsilon p(p-1)}{2}y(t)+\frac{C_f}{\varepsilon^{p-1}}+\frac{C_g (p-1)}{\varepsilon ^{\frac{p-2}{2}}}
+\int_{[-\tau,0]}\lambda(ds)x(t+s),
\]
where
\begin{equation}\label{definelambda}
\lambda(ds):=\nu(ds)\cdot\frac{1}{\varepsilon^{p-1}}+\eta(ds)\cdot \frac{p-1}{\varepsilon^{\frac{p-2}{2}}}.
\end{equation}
Hence
\begin{equation}\label{ineqy}
y(t)\leq e^{\beta_1t}y(0)+\int_0^te^{\beta_1(t-u)}\bigg(\beta_2+\beta_3+\int_{[-\tau,0]}\lambda(ds)x(u+s)\bigg)\,du,
\end{equation}
where
\[
\beta_1:=\frac{\varepsilon p(p-1)}{2},\quad \beta_2:=\frac{C_f}{\varepsilon^{p-1}},\quad \beta_3:=\frac{C_g (p-1)}{\varepsilon ^{\frac{p-2}{2}}}.
\]
Now since
\[
|X(t)|\leq |X(t)-D(X_t)|+|D(X_t)|,
\]
again by \eqref{Maoineqp},
\[
|X(t)|^p\leq (1+\varepsilon^{\frac{1}{p-1}})^{p-1}(\frac{1}{\varepsilon}|D(X_t)|^p+|X(t)-D(X_t)|^p),
\]
it follows that
\begin{align*}
x(t)&\leq (1+\varepsilon^{\frac{1}{p-1}})^{p-1}\bigg(\frac{1}{\varepsilon}\mathbb{E}[|D(X_t)|^p]+y(t)\bigg)\\
&\leq \frac{(1+\varepsilon^{\frac{1}{p-1}})^{p-1}}{\varepsilon}C_D+\frac{(1+\varepsilon^{\frac{1}{p-1}})^{p-1}}{\varepsilon}
\int_{[-\tau,0]}\mu(ds)x(t+s)
\\&\qquad\qquad\qquad\qquad\qquad+(1+\varepsilon^{\frac{1}{p-1}})^{p-1}y(t),
\end{align*}
Combining the above inequality with \eqref{ineqy}, we get
\begin{multline*}
x(t)\leq (1+\varepsilon^{\frac{1}{p-1}})^{p-1}e^{\beta_1t}y(0)+\beta_4C_D+\beta_4\int_{[-\tau,0]}\mu(ds)x(t+s)\\+(1+\varepsilon^{\frac{1}{p-1}})^{p-1}
\int_0^te^{\beta_1(t-u)}\bigg(\beta_2+\beta_3+\int_{[-\tau,0]}\lambda(ds)x(u+s)\bigg)\,du,
\end{multline*}
where $\beta_4:= (1+\varepsilon^{1/(p-1)})^{p-1}/\varepsilon$. Let $x_e(t)=e^{-\beta_1t}x(t)$ for $t\geq -\tau$. Since $e^{-\beta_1 t}\leq 1$ for $t\geq 0$, then
\begin{align*}
x_e(t)&\leq (1+\varepsilon^{\frac{1}{p-1}})^{p-1}y(0)+\beta_4C_D e^{-\beta_1t}+\beta_5(1-e^{-\beta_1t})
\\&\qquad\qquad\qquad
+\beta_4\int_{[-\tau,0]}\mu(ds)e^{-\beta_1t}x(t+s)
\\&\qquad\qquad\qquad\qquad\qquad
+(1+\varepsilon^{\frac{1}{p-1}})^{p-1}
\int_0^te^{-\beta_1u}\int_{[-\tau,0]}\lambda(ds)x(u+s)\,du
\\&\leq \bigg[(1+\varepsilon^{\frac{1}{p-1}})^{p-1}y(0)+\beta_4C_D+\beta_5\bigg]+\beta_4\int_{[-\tau,0]}
e^{\beta_1s}\mu(ds)x_e(t+s)
\\&\qquad\qquad\qquad\qquad\qquad\qquad\qquad+\int_0^t\int_{[-\tau,0]}e^{\beta_1s}\lambda (ds)x_e(u+s)\,du,
\end{align*}
where
\[
\beta_5:=\frac{1}{\beta_1}(1+\varepsilon^{\frac{1}{p-1}})^{p-1}(\beta_2+\beta_3).
\]
Let $\beta_6:=(1+\varepsilon^{\frac{1}{p-1}})^{p-1}y(0)+\beta_4C_D+\beta_5$, $\mu_e(ds):=e^{\beta_1s}\mu(ds)$ and
$\lambda_e(ds):=e^{\beta_1s}\lambda(ds)$, thus
\[
x_e(t)\leq \beta_6+\beta_4\int_{[-\tau,0]}\mu_e(ds)x_e(t+s)+\int_0^t\int_{[-\tau,0]}\lambda_e(ds)x_e(u+s)\,du.
\]
Now let $\mu(E)=\lambda(E)=0$ for $E \subset (-\infty, -\tau)$, so $\mu_e(E)=\lambda_e(E)=0$ for $E \subset (-\infty, -\tau)$.
Define $\mu^+_e(E):=\mu_e(-E)$ and $\lambda^+_e(E):=\lambda_e(-E)$ for $E \subset [0,\infty)$. Hence
\begin{align*}
\int_{[-\tau,0]}\mu_e(ds)x_e(t+s)&=\int_{(-\infty,0]}\mu_e(ds)x_e(t+s)
\\&= \int_{[0,\infty)}\mu^+_e(ds)x_e(t-s)
\\&= \int_{[0,t]}\mu^+_e(ds)x_e(t-s)+\int_{(t,\infty)}\mu^+_e(ds)x_e(t-s)
\\&= \int_{[0,t]}\mu^+_e(ds)x_e(t-s)+\int_{(t, t+\tau]}\mu^+_e(ds)\psi_e(t-s),
\end{align*}
where $\psi_e(t):=e^{-\beta_1 t}|\psi(t)|^p$ and $\psi$ is the initial condition for $X$ on $[-\tau,0]$. Similarly,
\[
\int_{[-\tau,0]}\lambda_e(ds)x_e(u+s)=\int_{[0,t]}\lambda^+_e(ds)x_e(u-s)+\int_{(t, t+\tau]}\lambda^+_e(ds)\psi_e(u-s).
\]
Consequently,
\begin{multline}\label{ineqxe}
x_e(t)\leq \beta_6+\beta_4\int_{[0,t]}\mu^+_e(ds)x_e(t-s)+\beta_4\int_{(t,t+\tau]}\mu^+_e(ds)\psi_e(t-s)
\\ + \int_0^t\int_{[0,u]}\lambda^+_e(ds)x_e(u-s)\,du+\int_0^t\int_{(u,u+\tau]}\lambda^+_e(ds)\psi_e(u-s)\,du.
\end{multline}
Let $\Lambda^+_e(t):=\int_{[0,t]}\lambda^+_e(ds)$. By Fubini's theorem and the integration-by-parts formula,
\begin{align}\label{lam1}
\int_{0}^t\int_{[0,u]}\lambda^+_e(ds)x_e(u-s)\,du
&= \int_{s=0}^t\lambda^+_e(ds)\int_{u=s}^t x_e(u-s)\,du
\\ \nonumber &= \int_{s=0}^t\lambda^+_e(ds)\int_{v=0}^{t-s}x_e(v)\,dv
\\ \nonumber &= \int_0^{t-s}x^+_e(v)\,dv\cdot \Lambda^+_e(s)\bigg|_{s=0}^t+\int_0^t\Lambda^+_e(s)x_e(t-s)\,ds\\
\nonumber &= \int_0^t\Lambda^+_e(s)x_e(t-s)\,ds.
\end{align}
Also
\begin{align}\nonumber
&\int_0^t\int_{(u,u+\tau]}\lambda^+_e(ds)\psi_e(u-s)\,du
\\ \nonumber &=\int_{[0,t+\tau]}\lambda^+_e(ds)\int_{(s-\tau)\vee 0}^{s\wedge t}\psi_e(u-s)\,du
\\ \nonumber &=\int_{[0,t]}\lambda^+_e(ds)\int_{(s-\tau)\vee 0}^{s\wedge t}\psi_e(u-s)\,du+\int_{(t,t+\tau]}\lambda^+_e(ds)\int_{(s-\tau)\vee 0}^{s\wedge t}\psi_e(u-s)\,du
\\ \label{lam2} &=\int_{[0,t]}\lambda^+_e(ds)\int_{(s-\tau)\vee 0}^s\psi_e(u-s)\,du+\int_{(t,t+\tau]}\lambda^+_e(ds)\int_{(s-\tau)\vee 0}^{t}\psi_e(u-s)\,du.
\end{align}
Now, if $t\geq \tau$, the second integral in \eqref{lam2} is zero; if $0\leq t< \tau$, then
\begin{align}\label{secondint}
\int_{(t,t+\tau]}\lambda^+_e(ds)\int_{(s-\tau)\vee 0}^{t}\psi_e(u-s)\,du
&= \int_{(t,\tau]}\lambda^+_e(ds)\int_{(s-\tau)\vee 0}^{t}\psi_e(u-s)\,du
\\ \nonumber &= \int_{(t,\tau]}\lambda^+_e(ds)\int_0^{t}\psi_e(u-s)\,du
\\ \nonumber &= \int_{(t,\tau]}\lambda^+_e(ds)\int_{-s}^{t-s}\psi_e(v)\,dv
\\ \nonumber &\leq \int_{(t,\tau]}\lambda^+_e(ds)\tau ||\psi_e||_{\text{sup}}
\\ \nonumber &\leq \tau ||\psi_e||_{\text{sup}}\int_{[0,\tau]}\lambda^+_e(ds).
\end{align}
For the first integral in \eqref{lam2},
\begin{align}\label{firstint}
\int_{[0,t]}\lambda^+_e(ds)\int_{(s-\tau)\vee 0}^s\psi_e(u-s)\,du
&= \int_{[0,\tau]}\lambda^+_e(ds)\int_{(s-\tau)\vee 0}^s\psi_e(u-s)\,du
\\ \nonumber &= \int_{[0,\tau]}\lambda^+_e(ds)\int_{0}^s\psi_e(u-s)\,du
\\ \nonumber &= \int_{[0,\tau]}\lambda^+_e(ds)\int_{-s}^0\psi_e(v)\,dv
\\ \nonumber &\leq \int_{[0,\tau]}\lambda^+_e(ds)\tau||\psi_e||_{\text{sup}}
\\ \nonumber &\leq \tau ||\psi_e||_{\text{sup}}\int_{[0,\tau]}\lambda^+_e(ds).
\end{align}
Inserting \eqref{secondint} and \eqref{firstint} into \eqref{lam2}, we have
\begin{equation}\label{ineqlam2}
\int_0^t\int_{(u,u+\tau]}\lambda^+_e(ds)\psi_e(u-s)\,du\leq 2\tau ||\psi_e||_{\text{sup}}\int_{[0,\tau]}\lambda^+_e(ds).
\end{equation}
Moreover, if $t\geq \tau$, then
\[
\beta_4\int_{(t,t+\tau]}\mu^+_e(ds)\psi_e(t-s)=0;
\]
if $0\leq t<\tau$, then
\begin{align}\label{mu}
\beta_4\int_{(t,t+\tau]}\mu^+_e(ds)\psi_e(t-s)
&= \beta_4\int_{(t,\tau]}\mu^+_e(ds)\psi_e(t-s)
\\ \nonumber &\leq \beta_4\int_{(t,\tau]}\mu^+_e(ds)||\psi_e||_{\text{sup}}
\\ \nonumber &\leq \beta_4 ||\psi_e||_{\text{sup}} \int_{[0,\tau]}\mu^+_e(ds).
\end{align}
Therefore combining \eqref{lam1}, \eqref{ineqlam2} and \eqref{mu} with \eqref{ineqxe}, we have
\begin{equation}
x_e(t)\leq \beta_7+\int_{[0,t]}\bigg(\beta_4\mu^+_e(ds)+\Lambda_e^+(s)\,ds\bigg)x_e(t-s),\quad t\geq 0.
\end{equation}
where
\[
\beta_7:=\beta_6 +\bigg(\beta_4 \int_{[0,\tau]}\mu^+_e(ds)+2\tau \int_{[0,\tau]}\lambda^+_e(ds)\bigg)||\psi_e||_{\text{sup}}.
\]
Choose small $\rho>0$ and define
\[
z(t):=\beta_7+\int_{[0,t]}\bigg(\beta_4\mu^+_e(ds)+\Lambda_e^+(s)\,ds+\rho\,ds\bigg)z(t-s),\quad t\geq 0.
\]
Then by Lemma \ref{lemagronwall}, we get $z(t)\geq x_e(t)$ for $t\geq 0$.

Next we determine the asymptotic behaviour of $z$. Note that the measure
\begin{equation}\label{definealpha}
\alpha(ds):=\beta_4\mu^+_e(ds)+\Lambda_e^+(s)\,ds+\rho\,ds
\end{equation}
has an absolutely continuous component. Moreover $\alpha$ is a positive measure. Also, we can find a number $\delta>0$ such that
$\int_{[0,\infty)} e^{-\delta s}\alpha(ds)=1$. Now, define the measure $\alpha_{\delta}\in M([0,\infty),\mathbb{R})$ by $\alpha_{\delta}(ds)=e^{-\delta s}\alpha(ds)$. Then $\alpha_{\delta}$ is a positive measure with a nontrivial absolutely continuous component such that $\alpha_{\delta}(\mathbb{R}^+)=1$. Also, we have that
\begin{align*}
\int_{[0,\infty)} s\alpha_{\delta}(ds)
&=\int_{[0,\infty)} se^{-\delta s}\alpha(ds)\\
&=\int_{[0,\infty)} se^{-\delta s}(\beta_4\mu^+_e(ds)+\Lambda_e^+(s)\,ds+\rho\,ds)\\
&=\beta_4\int_{[0,\tau]} se^{-\delta s}\mu^+_e(ds) + \int_{[0,\infty)} se^{-\delta s}\Lambda_e^+(s)\,ds+\rho\int_{[0,\infty)} se^{-\delta s}\,ds,
\end{align*}
since $\mu^+_e(E)=0$ for all $E\subset (\tau,\infty)$. Now, we note that because $\Lambda_e^+(t)\leq \Lambda_e^+(\infty)=\int_{[0,\tau]} \lambda_e^+(ds)<+\infty$ for all $t\geq 0$, the second integral on the righthand side is finite, and therefore we have that
$\int_{[0,\infty)} t\alpha_{\delta}(dt)<+\infty$. Next define $z_{\delta}(t):=e^{-\delta t} z(t)$ for $t\geq 0$
so that
\[
z_{\delta}(t)=\beta_7 e^{-\delta t}+\int_{[0,t]} \alpha_{\delta}(ds) z_{\delta}(t-s),\quad t\geq 0.
\]
Now, define $-\gamma$ to be the resolvent of $-\alpha_{\delta}$. Then, by the renewal theorem (see Gripenberg, Londen and Staffans~\cite[Theorem 7.4.1]{GrLoSt90}), the existence of $\gamma$ is guaranteed. Moreover, $\gamma$ is a positive measure and is of the form
\[
\gamma(dt)=\gamma_1(dt)+\gamma_1([0,t])\,dt
\]
where $\gamma_1\in M(\mathbb{R}^+;\mathbb{R})$ and $\gamma_1(\mathbb{R}^+)=1/\int_{\mathbb{R}^+} t\alpha_{\delta}(dt)$, which is finite. Since
$(-\gamma)+(-\alpha_{\delta})\ast(-\gamma)=-\alpha_{\delta}$, let $h(t):=\beta_7 e^{-\delta t}$, we have
\begin{align*}
z_\delta=h+\alpha_\delta\ast z_\delta
&=h+\gamma\ast z_\delta-\alpha_\delta\ast \gamma \ast z_\delta
\\&= h+\gamma\ast (z_\delta-\alpha_\delta\ast z_\delta)
\\&=h+\gamma\ast h,
\end{align*}
that is
\begin{align*}
z_\delta(t)&=\beta_7 e^{-\delta t}+\beta_7 \int_{[0,t]}\gamma(ds)e^{-\delta(t-s)}
\\&= \beta_7 e^{-\delta t}+\beta_7 \int_{[0,t]}\bigg(\gamma_1(ds)+\gamma_1([0,s])ds\bigg)e^{-\delta(t-s)}.
\end{align*}
Thus
\[
\limsup_{t\to\infty}\frac{x_e(t)}{e^{\delta t}}\leq \limsup_{t\to\infty}\frac{z(t)}{e^{\delta t}}
=\limsup_{t\to\infty}z_\delta(t)\leq \frac{\beta_7}{\int_{\mathbb{R^+}}t\alpha_\delta(dt)}+\frac{\beta_7}{\delta
\int_{\mathbb{R^+}}t\alpha_\delta(dt)}.
\]
Hence there exists $C>0$ such that $x_e(t)\leq C e^{\delta t}$ for $t\geq 0$. Therefore $\mathbb{E}[|X(t)|^p]=x(t)=e^{\beta_1 t}x_e(t)\leq Ce^{(\delta +\beta_1)t}$ for $t\geq 0$, which implies
\[
\limsup_{t\to\infty}\frac{1}{t}\log{\mathbb{E}[|X(t)|^p]}\leq \delta+\beta_1.
\]
Now in \eqref{definealpha}, let $\rho\to 0$, then $\delta \to \delta_\ast$, where
\begin{equation}\label{mass1}
\int_{[0,\infty)} e^{-\delta_\ast s}\alpha(ds)= \int_{[0,\infty]}e^{-\delta_\ast s}\bigg(\mu^+_e(ds)+\Lambda^+_e(s)\,ds\bigg)=1.
\end{equation}
Note
\begin{multline*}
\int_0^\infty e^{-\delta_\ast s}\mu^+_e(ds)=\int_0^\tau e^{-\delta_\ast s}\mu^+_e(ds)=\int^0_{-\tau} e^{-\delta_\ast s}\mu_e(ds)
=\int ^0_{-\tau} e^{(-\delta_\ast +\beta_1)s}\mu(ds),
\end{multline*}
and
\begin{align*}
\int_0^\infty e^{-\delta_\ast s} \Lambda^+_e(s)\,ds
&= \int_0^\tau e^{-\delta_\ast s}\Lambda^+_e(s)\,ds+\int_\tau^\infty e^{-\delta_\ast s}\int_{[0,\tau]}\lambda^+_e(du)\,ds
\\&= \int_0^\tau e^{-\delta_\ast s}\int_{[0,s]}\lambda^+_e(du)\,ds+\frac{e^{-\delta_\ast \tau}}{\delta_\ast}\int_{[0,\tau]}\lambda^+_e(du)
\\&= \int_0^\tau e^{-\delta_\ast s}\int_{[-s,0]}\lambda_e(du)\,ds+\frac{e^{-\delta_\ast \tau}}{\delta_\ast}\int_{[-\tau,0]}\lambda_e(du)
\\&= \int_0^\tau e^{-\delta_\ast s}\int_{[-s,0]}e^{\beta_1 u}\lambda(du)\,ds+\frac{e^{-\delta_\ast \tau}}{\delta_\ast}\int_{[-\tau,0]}e^{\beta_1 u}\lambda(du).
\end{align*}
where $\lambda$ is defined in \eqref{definelambda}. Replace $\delta_\ast$ by $\delta$, we get the desired result.

Finally, we show that \eqref{finite4thmome} holds for $t\geq 0$. By H\"{o}lder's inequality, we get
\[
\mathbb{E}\bigg[\int_0^t|Y(s)|^{2p-2}||g(X_s)||^2\bigg]\,ds\leq \int_0^t\mathbb{E}[|Y(s)|^{4p-4}]^\frac{1}{2}\mathbb{E}[||g(X_s)||^4]^\frac{1}{2}\,ds.
\]
Given \eqref{g}, by Lemma \ref{lemafinitemo}, let $\varepsilon=1$ in \eqref{Maoineqp}, there exist positive real numbers $K_1$ and $K_2$ such that
\begin{align*}
\mathbb{E}[||g(X_s)||^4]&\leq \mathbb{E}\bigg[\bigg(C_g+\int_{[-\tau,0]}\eta(du)|X(s+u)|\bigg)^4\bigg]
\\&\leq 8\mathbb{E}\bigg[C^4_g+\bigg(\int_{[-\tau,0]}\eta(du)|X(s+u)|\bigg)^4\bigg]
\\&\leq 8C^4_g+8\bigg(\int_{[-\tau,0]}\eta(du)\bigg)^3\bigg(\int_{[-\tau,0]}\eta(du)\mathbb{E}[|X(s+u)|^4]\bigg)
\\&\leq 8C^4_g+8\bigg(\int_{[-\tau,0]}\eta(du)\bigg)^4K_1e^{K_2s}.
\end{align*}
There also exist positive real numbers $K_3$ and $K_4$ such that
\begin{align*}
\mathbb{E}[|Y(s)|^{4p-4}]&=\mathbb{E}[|X(s)-D(X_s)|^{4p-4}]
\\&\leq 2^{4p-5}\bigg(\mathbb{E}[|X(s)|^{4p-4}]+\mathbb{E}[|D(X_s)|^{4p-4}]\bigg)
\\&\leq 2^{4p-5}\bigg(K_3e^{K_4 s}+\mathbb{E}[|D(X_s)|^{4p-4}]\bigg).
\end{align*}
Apply the same analysis to $\mathbb{E}[|D(X_s)|^{4p-4}]$ as $\mathbb{E}[||g(X_s)||^4]$ using \eqref{D}, it is easy to see that
\[
\int_0^t\mathbb{E}[|Y(s)|^{4p-4}]^\frac{1}{2}\mathbb{E}[||g(X_s)||^4]^\frac{1}{2}\,ds<\infty.
\]
Hence \eqref{finite4thmome} holds.

\subsection{Proof of Theorem \ref{thmexpestas}}
Again, let $Y(t):=X(t)-D(X_t)$ for any $t\geq 0$. For any $n\leq t\leq n+1$, we have
\[
|Y(t)|^2\leq \frac{(1+\varepsilon)^2}{\varepsilon^2}|Y(n)|^2+\frac{(1+\varepsilon)^2}{\varepsilon}\left(\int_n^t|f(X_s)|\,ds\right)^2
+(1+\varepsilon)\left|\int_n^tg(X_s)dB(s)\right|^2.
\]
Therefore by \eqref{f1} we have 
\begin{align*}
&\sup_{n\leq t\leq n+1} |Y(t)|^2
\\&\leq \frac{(1+\varepsilon)^2}{\varepsilon^2}|Y(n)|^2+
\frac{(1+\varepsilon)^2}{\varepsilon}\left(\int_n^{n+1} |f(X_s)|\,ds\right)^2\\
&\qquad\qquad+(1+\varepsilon)\sup_{n\leq t\leq n+1} \bigg|\int_n^tg(X_s)dB(s)\bigg|^2\\
&\leq \frac{(1+\varepsilon)^2}{\varepsilon^2}|Y(n)|^2+\frac{(1+\varepsilon)^2}{\varepsilon}\int_n^{n+1} |f(X_s)|^2\,ds
\\
&\qquad\qquad+(1+\varepsilon)\sup_{n\leq t\leq n+1} \bigg|\int_n^tg(X_s)dB(s)\bigg|^2\\
&\leq \frac{(1+\varepsilon)^2}{\varepsilon^2}|Y(n)|^2+\frac{(1+\varepsilon)^2}{\varepsilon}\int_n^{n+1} \left(C_f+\int_{[-\tau,0]} \nu(du)|X(s+u)|\right)^2\,ds\\
&\qquad\qquad+(1+\varepsilon)\sup_{n\leq t\leq n+1} \bigg|\int_n^tg(X_s)dB(s)\bigg|^2\\
&\leq \frac{(1+\varepsilon)^2}{\varepsilon^2}|Y(n)|^2+\frac{(1+\varepsilon)^3}{\varepsilon^2}C_f^2+\frac{(1+\varepsilon)^3}{\varepsilon}\int_n^{n+1} \bigg(\int_{[-\tau,0]} \nu(du)|X(s+u)|\bigg)^2\,ds\\
&\qquad\qquad+(1+\varepsilon)\sup_{n\leq t\leq n+1} \bigg|\int_n^tg(X_s)dB(s)\bigg|^2\\
&\leq \frac{(1+\varepsilon)^2}{\varepsilon^2}|Y(n)|^2+\frac{(1+\varepsilon)^3}{\varepsilon^2}C_f^2
+(1+\varepsilon)\sup_{n\leq t\leq n+1} \bigg|\int_n^tg(X_s)dB(s)\bigg|^2\\
&\qquad+\frac{(1+\varepsilon)^3}{\varepsilon}\bigg(\int_{[-\tau,0]} \nu(du)\bigg)\int_n^{n+1} \int_{[-\tau,0]}\nu(du)|X(s+u)|^2\,ds.
\end{align*}
Hence by the Burkholder--Davis--Gundy inequality, with $x(t):=\mathbb{E}[|X(t)|^2]$, we have 
\begin{align*}
&\mathbb{E}\left[\sup_{n\leq t\leq n+1} \bigg|\int_n^tg(X_s)dB(s)\bigg|^2\right]
\\&\leq 4\mathbb{E}\left[\int_n^{n+1} \|g(X_s)\|^2\,ds\right]\\
&\leq 4\mathbb{E}\left[\int_n^{n+1} \left(C_g+\int_{[-\tau,0]} \eta(du)|X(s+u)|\right)^2\,ds\right]\\
&\leq 4\frac{1+\varepsilon}{\varepsilon}C_g^2+4(1+\varepsilon)\bigg(\int_{[-\tau,0]}\eta(du)\bigg)
\mathbb{E}\bigg[\int_n^{n+1}\int_{[-\tau,0]}\eta(du)|X(s+u)|^2\,ds\bigg]
\\&=4\frac{1+\varepsilon}{\varepsilon}C_g^2+4(1+\varepsilon)\bigg(\int_{[-\tau,0]}\eta(du)\bigg)
\int_n^{n+1}\int_{[-\tau,0]}\eta(du)x(s+u)\,ds.
\end{align*}
Therefore we deduce that 
\begin{align}
\mathbb{E}\left[\sup_{n\leq t\leq n+1} |Y(t)|^2\right]
&\leq \frac{(1+\varepsilon)^2}{\varepsilon^2}\mathbb{E}|Y(n)|^2+\frac{(1+\varepsilon)^3}{\varepsilon^2}C_f^2
+
4\frac{(1+\varepsilon)^2}{\varepsilon}C_g^2\nonumber\\
&\quad+4(1+\varepsilon)^2\bigg(\int_{[-\tau,0]}\eta(du)\bigg)\int_n^{n+1}\int_{[-\tau,0]}\eta(du)x(s+u)\,ds\nonumber\\
 \label{eq.estexpsupY}
&\qquad+\frac{(1+\varepsilon)^3}{\varepsilon}\bigg(\int_{[-\tau,0]} \nu(du)\bigg)\int_n^{n+1} \int_{[-\tau,0]}\nu(du)x(s+u)\,ds.
\end{align}
$x(t)$ is bounded above by an exponential. We now show that $\mathbb{E}[Y(t)^2]$ can also be bounded by an exponential.  
Since $Y(t)=X(t)-D(X_t)$, and $D$ obeys \eqref{D1} we have 
\begin{align*}
|Y(t)|^2 &\leq \frac{1+\varepsilon}{\varepsilon}|X(t)|^2+(1+\varepsilon)\left(C_D+\int_{[-\tau,0]} \mu(ds)|X(t+s)|\right)^2\\
&\leq \frac{1+\varepsilon}{\varepsilon}|X(t)|^2+\frac{(1+\varepsilon)^2}{\varepsilon} C_D^2
\\
&\qquad\qquad +  (1+\varepsilon)^2\int_{[-\tau,0]} \mu(ds) \int_{[-\tau,0]} \mu(ds)|X(t+s)|^2,
\end{align*}
so
\begin{equation}  \label{eq.estexpY}
\mathbb{E}[|Y(t)|^2]
\leq \frac{1+\varepsilon}{\varepsilon}x(t)+\frac{(1+\varepsilon)^2}{\varepsilon} C_D^2
+  (1+\varepsilon)^2\int_{[-\tau,0]} \mu(ds) \int_{[-\tau,0]} \mu(ds)x(t+s).
\end{equation}
By Theorem~\ref{thmexpest}, there exist $C_0>0$ and $\gamma>0$ such that $x(t)=\mathbb{E}[|X(t)|^2]\leq C_0e^{\gamma t}$ for all $t\geq -\tau$. 
Inserting this estimate in \eqref{eq.estexpY} shows that there is $C_1>0$ such that $\mathbb{E}[|Y(t)|^2]\leq C_1 e^{\gamma t}$ for all $t\geq 0$.
Using this estimate and $x(t)\leq C_0e^{\gamma t}$ for all $t\geq -\tau$ in \eqref{eq.estexpsupY}, we see that there is a $C_2>0$ such that 
\begin{equation*} 
\mathbb{E}\left[\sup_{n\leq t\leq n+1} |Y(t)|^2\right]\leq C_2 e^{\gamma n}, \quad \text{for all $n\geq 0$.}
\end{equation*}
By the Borel--Cantelli lemma, it then follows that  
\[
\limsup_{t\to\infty} \frac{1}{t}\log |Y(t)|\leq \frac{\gamma}{2},\quad\text{a.s.}
\]
Therefore, as $t\mapsto |Y(t)|$ is continuous almost surely, for every $\varepsilon>0$, there exists an almost surely finite random variable $C_3=C_3(\varepsilon)>0$ such that
\begin{equation}  \label{eq.estexpasY}
|Y(t)|\leq C_3(\varepsilon) e^{(\gamma/2+\varepsilon) t}, \quad \text{for all $t\geq 0$, a.s.} 
\end{equation}

We are finally in a position to estimate $|X(t)|$. By \eqref{D1} we have 
\[
|X(t)|\leq |Y(t)| + C_D + \int_{[-\tau,0]} \mu(ds)|X(t+s)|.
\]
Define $\mu^+(E):=\mu(-E)$ for all $E\subseteq (-\infty,0]$ where we have extended $\mu$ to all of $(-\infty,0]$ as before. Then
\begin{align*}
|X(t)|&\leq |Y(t)| + C_D + \int_{[0,\infty)} \mu^+(ds)|X(t-s)|\\
&\leq |Y(t)| + C_D + \int_{[0,t]} \mu^+(ds)|X(t-s)|+\int_{(t,t+\tau]} \mu^+(ds)|\phi(t-s)|,
\end{align*}
where we have once again extended $X$ to be zero on $(-\infty,-\tau)$. Clearly, the last term is bounded, so by \eqref{eq.estexpasY}
there exists an almost surely finite random variable $C_4=C_4(\varepsilon)>0$ such that 
\begin{equation} \label{eq.xvolterraeqnas}
|X(t)|
\leq C_4(\varepsilon) e^{(\gamma/2+\varepsilon) t}+\int_{[0,t]} \mu^+(ds)|X(t-s)|,\quad t\geq 0.
\end{equation}

Now let $\rho>0$ and consider $Z$ which is defined as
\begin{equation*} 
Z(t)=C_4(\varepsilon)e^{(\gamma/2+\varepsilon) t}+\int_{[0,t]} \mu^{+}_\rho(ds)Z(t-s),\quad t\geq 0.
\end{equation*}
where $\mu^{+}_\rho(ds):=\mu^+(ds)+ \rho e^{-s} \,ds$.  
Then by Lemma \ref{lemagronwall}, $Z(t)\geq |X(t)|$ for $t\geq 0$ a.s. Clearly $Z$ is positive also. 

We consider first the case when $\mu^+(\mathbb{R}^+)\geq 1$. Let $\theta>0$ and define 
\[
\alpha_{\theta,\rho}(ds):=e^{-\theta s} \mu^{+}_\rho(ds).
\]
Note for all $\theta$ and $\rho$ that $\alpha_\theta$ has an absolutely continuous component. Also we note that 
%\begin{itemize}
%\item[(i)] If $\mu^+(\mathbb{R}^+)<1$ and $0<\rho<1-\mu^+(\mathbb{R}^+)$, then $\mu^{+}_\rho(\mathbb{R}^+)<1$;
Since $\rho>0$ if $\mu^+(\mathbb{R}^+)=1$ then $\mu^{+}_\rho(\mathbb{R}^+)>1$; obviously 
if $\mu^+(\mathbb{R}^+)>1$, then $\mu^{+}_\rho(\mathbb{R}^+)>1$. Therefore, in each case there exists $\theta>0$ such that $\int_{[0,\infty)}\alpha_{\theta,\rho}(ds)=1$. 
%In the case when $\mu^+(\mathbb{R}^+)\geq 1$, we have $\theta>0$, while $\mu^+(\mathbb{R}^+)< 1$ implies $\theta<0$.
Next we show that $\int_{[0,\infty)}s\alpha_{\theta,\rho}(ds)<+\infty$. This follows from  
\begin{equation*}
\int_{[0,\infty)}s\alpha_{\theta,\rho}(ds)
%\int_{[0,\infty)}se^{-\theta s}\bigl(\mu^+(ds)+ \rho e^{-s} \,ds\bigr) \\
=
\int_{[0,\tau]}se^{-\theta s}\mu^+(ds)+ \rho\int_{[0,\infty)} se^{-(\theta+1) s}  \,ds<\infty.
\end{equation*}
Let $Z_\theta(t):=Z(t)e^{-\theta t}$. Hence
\[
Z_\theta(t)=C_4(\varepsilon) e^{(\frac{\gamma}{2}+\varepsilon-\theta)t}+\int_{[0,t]}\alpha_{\theta,\rho}(ds)Z_\theta(t-s), \quad t\geq 0.
\]
Applying the same argument as in the proof of Theorem \ref{thmexpest}, there exists a measure $\gamma$ such that $-\gamma$ is the resolvent
of $-\alpha_{\theta,\rho}$. Moreover, $\gamma(dt)=\gamma_1(dt)+\gamma_1([0,t])dt$, where $\gamma_1\in M(\mathbb{R}^+;\mathbb{R})$ and $\gamma_1(\mathbb{R}^+)
=(\int_{\mathbb{R}^+}t\alpha_{\theta,\rho}(dt))^{-1}$. Hence
\begin{align}
Z_\theta(t)&=C_4(\varepsilon)e^{(\frac{\gamma}{2}+\varepsilon-\theta)t}+C_4(\varepsilon)\int_{[0,t]}\gamma(ds)e^{(\frac{\gamma}{2}+\varepsilon-\theta)
(t-s)}\nonumber\\
&=C_4(\varepsilon)e^{(\frac{\gamma}{2}+\varepsilon-\theta)t}+C_4(\varepsilon)\int_{[0,t]}\bigg(\gamma_1(ds)+\gamma_1([0,s])\,ds\bigg)
e^{(\frac{\gamma}{2}+\varepsilon-\theta)(t-s)}\nonumber\\
&=C_4(\varepsilon)e^{(\frac{\gamma}{2}+\varepsilon-\theta)t}
+C_4(\varepsilon)\int_{[0,t]} e^{(\frac{\gamma}{2}+\varepsilon-\theta)(t-s)}\gamma_1(ds)\nonumber\\
\label{eq.zintestdiffth}
&\qquad+C_4(\varepsilon)\int_{[0,t]}e^{(\frac{\gamma}{2}+\varepsilon-\theta)(t-s)}\gamma_1([0,s])\,ds.
\end{align}
In the case when $\theta>\gamma/2+\varepsilon$, we have 
\[
\limsup_{t\to\infty} \frac{|X(t)|}{e^{\theta t}}\leq 
\lim_{t\to\infty} \frac{Z(t)}{e^{\theta t}}
=\lim_{t\to\infty} Z_\theta(t)=\frac{C_4(\varepsilon)}{(\theta-\gamma/2-\varepsilon)\int_{[0,\infty)} t\alpha_{\theta,\rho}(dt)}, 
\quad \text{a.s.}
\]
If $\theta=\gamma/2+\varepsilon$, then 
\[
Z_\theta(t)
=C_4(\varepsilon)+C_4(\varepsilon)\int_{[0,t]} \gamma_1(ds)+C_4(\varepsilon)\int_{[0,t]}\gamma_1([0,s])\,ds,
\]
so
\[
\limsup_{t\to\infty} \frac{|X(t)|}{te^{\theta t}}\leq
\lim_{t\to\infty} \frac{Z(t)}{te^{\theta t}}
=\lim_{t\to\infty} \frac{Z_\theta(t)}{t}
=\frac{C_4(\varepsilon)}{\int_{[0,\infty)} t\alpha_{\theta,\rho}(dt)}, \quad\text{a.s.}
\]
In the case when $\theta<\gamma/2+\varepsilon$, we have 
\begin{multline*}
\frac{Z_\theta(t)}{e^{(\frac{\gamma}{2}+\varepsilon-\theta)t}}=C_4(\varepsilon)
+C_4(\varepsilon)\int_{[0,t]} e^{-(\frac{\gamma}{2}+\varepsilon-\theta)s}\gamma_1(ds)
%\label{eq.zintestdiffth}
\\+C_4(\varepsilon)\int_{[0,t]}e^{-(\frac{\gamma}{2}+\varepsilon-\theta)s}\gamma_1([0,s])\,ds.
\end{multline*}
Therefore
\[
\lim_{t\to\infty} \frac{Z_\theta(t)}{e^{(\frac{\gamma}{2}+\varepsilon-\theta)t}}=C_4(\varepsilon)
+C_4(\varepsilon)\int_{[0,\infty)
} e^{-(\frac{\gamma}{2}+\varepsilon-\theta)s}\left(\gamma_1(ds)+\gamma_1([0,s])\right)\,ds.
\]
Hence
\begin{align*}
\limsup_{t\to\infty}\frac{|X(t)|}{e^{(\gamma/2+\varepsilon)t}}
&\leq \lim_{t\to\infty}\frac{Z(t)}{e^{(\gamma/2+\varepsilon)t}}=\lim_{t\to\infty} \frac{Z_\theta(t)}{{e^{(\gamma/2+\varepsilon-\theta)t}}}\\
&= C_4(\varepsilon)
+C_4(\varepsilon)\int_{[0,\infty)
} e^{-(\frac{\gamma}{2}+\varepsilon-\theta)s}\left(\gamma_1(ds)+\gamma_1([0,s])\right)\,ds,\quad\text{a.s.}
\end{align*}
We have therefore shown that $\mu^+(\mathbb{R}^+)\geq 1$ implies that 
\[
\limsup_{t\to\infty} \frac{1}{t}\log |X(t)|\leq \max(\gamma/2+\varepsilon,\theta(\rho)), \quad\text{a.s.}
\]
Since $\rho$ and $\varepsilon>0$ are arbitrary, we may send them both to zero through the rational numbers to get
\begin{equation} \label{eq.Xasest1mugt1}
\limsup_{t\to\infty} \frac{1}{t}\log |X(t)|\leq \max(\gamma/2,\theta^\ast), \quad\text{a.s.}
\end{equation}
where $\theta^\ast\geq 0$ is such that $\int_{[0,\tau]} e^{-\theta^\ast s} \mu^{+}(ds) =1$.
 
We now consider the case where $\mu^+(\mathbb{R}^+)<1$. Recall that  
\begin{equation*}
|X(t)|
\leq C_4(\varepsilon) e^{(\gamma/2+\varepsilon) t}+\int_{[0,t]} \mu^+(ds)|X(t-s)|,\quad t\geq 0.
\end{equation*}
Therefore
\begin{equation*}
|X(t)| e^{-(\gamma/2+\varepsilon) t}
\leq C_4(\varepsilon)+\int_{[0,t]} e^{-(\gamma/2+\varepsilon)s} \mu^+(ds) e^{-(\gamma/2+\varepsilon)(t-s)}|X(t-s)|,\quad t\geq 0.
\end{equation*}
Since $\gamma/2+\varepsilon>0$, we have that $\int_{[0,\infty)} e^{-(\gamma/2+\varepsilon)s} \mu^+(ds)<1$. Therefore, there 
exists a positive $\rho:=1-\int_{[0,\infty)} e^{-(\gamma/2+\varepsilon)s} \mu^+(ds)$. Define 
\[
\mu_\rho(ds):=e^{-(\gamma/2+\varepsilon)s} \mu^+(ds) + \rho e^{-s}\,ds.
\]
Then $\mu_\rho$ is a positive measure with an absolutely continuous component and we have
\[
\int_{[0,\infty)}\mu_\rho(ds)= \int_{[0,\infty)} e^{-(\gamma/2+\varepsilon)s} \mu^+(ds) + \int_0^\infty \rho e^{-s}\,ds
%=\int_{[0,\infty)} e^{-(\gamma/2+\varepsilon)s} \mu^+(ds)+\rho
=1,
\]
by the definition of $\rho$. Also 
\[
\int_{[0,\infty)} s\mu_\rho(ds)= \int_{[0,\infty)} se^{-(\gamma/2+\varepsilon)s} \mu^+(ds) + \int_0^\infty \rho se^{-s}\,ds<\infty.
\]
Now let $Z$ be the solution of 
\[
Z(t)=C_4(\varepsilon)+\int_{[0,t]} \mu_\rho(ds)Z(t-s),\quad t\geq 0.
\]
Clearly $Z$ is positive. Moreover $|X(t)| e^{-(\gamma/2+\varepsilon) t}\leq Z(t)$ for all $t\geq 0$ a.s. Let $-\gamma_1$ be the 
resolvent of $-\mu_\rho$. Then 
\[
Z(t)=C_4+C_4\int_{[0,t]} \gamma_1(ds), \quad t\geq 0.
\] 
By applying the same argument used in the proof of Theorem~\ref{thmexpest} we have that 
\[
\lim_{t\to\infty} \frac{Z(t)}{t}=C_4\gamma_1(\mathbb{R^+})=\frac{C_4}{\int_0^\infty t\mu_\rho(dt)}.
\]
Therefore 
\[
\limsup_{t\to\infty} \frac{|X(t)|}{te^{(\gamma/2+\varepsilon)t}}
\leq 
\lim_{t\to\infty} \frac{Z(t)}{t}=\frac{C_4}{\int_0^\infty t\mu_\rho(dt)}, \quad \text{a.s.},
\]
and so
\[
\limsup_{t\to\infty} \frac{1}{t}\log|X(t)|\leq \gamma/2+\varepsilon, \quad \text{a.s.}
\]
and so letting $\varepsilon\to 0^+$ through the rational numbers gives
\begin{equation} \label{eq.Xasest1mult1}
\limsup_{t\to\infty} \frac{1}{t}\log |X(t)|\leq \gamma/2, \quad\text{a.s.}
\end{equation}

\subsection{Proof of Proposition~\ref{prop.nonexiststochadd}}
Let $\Omega_1$ be an almost sure event such that $t\mapsto
B(t,\omega)$ is nowhere differentiable on $(0,\infty)$. Let $T>0$.
Suppose that $X=\{X(t):-\tau\leq t\leq T\}$ is a solution of \eqref{eq.maineqdiffnonexist}, \eqref{eq.icnonexist}.
Then $X$ is $(\mathcal{F}(t))_{t\geq 0}$--adapted and is such that $t\mapsto X(t,\omega)$ is continuous on $[-\tau,T]$ for all $\omega\in \Omega_2$, where $\Omega_2$ is an almost sure event.
Define $C_T=\{\omega:X(\cdot,\omega) \text{ obeys \eqref{eq.maineqintnonexist}}\}$ and
\[
A_T=C_T\cap\Omega_1\cap \Omega_2,
\]
Thus $\mathbb{P}[C_T]>0$ and so $\mathbb{P}[A_T]>0$. Hence for each $\omega\in A_T$,
we have for all $t\in[0,T]$
\[
\int_{-\tau}^0 w(s)h(X(t+s,\omega))\,ds = \int_{-\tau}^0
w(s)h(\psi(s))\,ds + \int_0^t f(X_s(\omega))\,ds + \sigma
B(t,\omega),
\]
so
\begin{equation} \label{eq.beqf}
\sigma B(t,\omega)= F(t,\omega), \quad t\in[0,T],
\end{equation}
where we have defined
\[
F(t,\omega):= \int_{-\tau}^0 w(s)h(X(t+s,\omega))\,ds
-\int_{-\tau}^0 w(s)h(\psi(s))\,ds - \int_0^t f(X_s(\omega))\,ds.
\]
It is not difficult to show that the righthand side of \eqref{eq.beqf} viz., $t\mapsto F(t,\omega)$ is
differentiable on $[0,T]$ for each $\omega\in A_T$, while the lefthand
side of \eqref{eq.beqf} is not differentiable anywhere in $[0,T]$ for each $\omega\in A_T$. This contradiction means that $\mathbb{P}[A_T]=0$;
hence with probability zero there are no sample paths of $X$ which satisfy \eqref{eq.maineqdiffnonexist}, \eqref{eq.icnonexist}.

\subsection{Proof of Proposition~\ref{prop.sfdemaxtypenonexist}}
Suppose $X$ is a solution on $[-\tau,T]$. Then with $A:=\psi(0)+\kappa \max_{s\in[-\tau,0]} |\psi(s)|$
\[
X(t)+\kappa \max_{s\in[t-\tau,t]} |X(s)|= A + \int_0^t g(X_s)\,dB(s), \quad t\in [0,T], \quad \text{a.s.}
\]
Clearly $X(t)+\kappa \max_{s\in[t-\tau,t]} |X(s)|\geq -|X(t)|+\kappa|X(t)|=(\kappa-1)|X(t)|\geq 0$. Therefore
\begin{equation} \label{eq.Mboundedmaxtypenonexist}
M(t):=
\int_0^t -g(X_s)\,dB(s) \leq A, \quad t\in [0,T], \quad \text{a.s.}
\end{equation}
Note that $A\geq 0$. Clearly $M$ is a local martingale with $\langle M \rangle (t)=\int_0^t g^2(X_s)\,ds \geq \delta t$ by \eqref{eq.nondegdiffsfdemaxtype}. By the martingale time change theorem, there exists a standard Brownian motion $\tilde{B}$ such that $M(t)=\tilde{B}(\langle M \rangle(t))$ for $t\in[0,T]$. Therefore by \eqref{eq.Mboundedmaxtypenonexist} we have
\[
\max_{0\leq u\leq T} \tilde{B}(\langle M\rangle(u))\leq A, \quad \text{a.s.}
\]
Since $\langle M \rangle (T) \geq \delta T$ and $t\mapsto \langle M \rangle(t)$ is increasing on $[0,T]$ we have
\[
\max_{0\leq s\leq \delta T} \tilde{B}(s) \leq \max_{0\leq u\leq T} \tilde{B}(\langle M\rangle(u))\leq A, \quad \text{a.s.},
\]
which is false, because $\tilde{B}$ is a standard Brownian motion $\delta T>0$ and $A\geq 0$ is finite, recalling that
$|W(\delta T)|$ and $\max_{0\leq s\leq \delta T} W(s)$ have the same distribution for any standard Brownian motion $W$.
Hence there is no process $X$ which is a solution on $[-\tau,T]$.

%%%%%%%%%%%%%%%%%%%%%%%%%%%%%%%%%%%%%%%%%%%%%%%%%%%%%%%%%%%%%%%%%%%%%%%%%%%%%%%%%%%%%%%%%%%%%%%%%%%%%%%%%%%%%%%%%%%%%%%%%%%%%%%%%%%%%%%%%%%%%%%%%%%

%%%%%%%%%%%%%%%%%%%%%%%%%%%%%%%%%%%%%%%%%%%%%%%%%%%%%%%%%%%%%%%%%%%%%%%%%%%%%%%%%%%%%%%%%%%%%%%%%%%%%%%%%%%%%%%%%%%%%%%%%%%%%%%%%%%%%%%%%%%%%%%%%%%%%%

%Neutral delay differential equations have been used to describe
%various processes in physics and engineering sciences \cite{HaleLun93}, \cite{Stepan:1989}. For
%example, transmission lines involving nonlinear boundary conditions
%\cite{Hale:77}, cell growth dynamics \cite{BakBochPaul:1998}, propagating pulses in cardiac tissue
%\cite{CourGlassKeen:1993} and drillstring vibrations \cite{BalJansMcClint:2003} have been described by means of
%neutral delay differential equations.

\end{document}